\documentclass[12pt]{amsart}
\usepackage{latexsym, amssymb, amscd, rotating}
\renewcommand\a{\alpha}
\renewcommand\b{\beta}
\newcommand\g{\gamma}
\renewcommand\d{\delta}
\newcommand\la{\lambda}
\newcommand\z{\zeta}
\newcommand\e{\eta}
\renewcommand\th{\theta}

\newcommand\m{\mu}

\newcommand\s{\sigma}
\newcommand\x{\chi}
\newcommand\f{\phi}
\newcommand\vf{\varphi}

\renewcommand\t{\tau}
\renewcommand\r{\rho}

\newcommand\Om{\Omega}
\newcommand\w{\omega}

\newcommand\vD{\varDelta}

\newcommand\F{\Phi}

\newcommand\vL{\varLambda}

\newcommand{\vT}{\varTheta}
\newcommand\vG{\varGamma}
\newcommand\ve{\varepsilon}

\newcommand{\QQ}{\mathbb Q}

\newcommand{\ZZ}{\mathbb Z}
\newcommand{\NN}{\mathbb N}
\newcommand{\CC}{\mathbb C}
\newcommand{\GG}{\mathbb G}

\newcommand\Fq{{\mathbf F}_q}

\newcommand\BQ{\mathbf Q}
\newcommand\BP{\mathbf P}

\newcommand\Bp{\boldsymbol p}
\newcommand\Bq{\boldsymbol q}

\newcommand\Bm{\mathbf m}
\newcommand\Bs{\mathbf s}
\newcommand\Bv{\mathbf v}
\newcommand\Bu{\mathbf u}

\newcommand\CE{\mathcal{E}}

\newcommand\CP{\mathcal{P}}

\newcommand\CF{\mathcal{F}}

\newcommand\CX{ \mathcal{X}}
\newcommand\CY{ \mathcal{Y}}

\newcommand\FS{\mathfrak S}

\newcommand\BLa{\boldsymbol\varLambda}

\newcommand\Ba{\boldsymbol\alpha}
\newcommand\Bb{\boldsymbol\beta}
\newcommand\Bg{\boldsymbol\gamma}

\newcommand\BOm{\boldsymbol\Omega}
\newcommand\iv{^{-1}}

\newcommand\wt{\widetilde}
\newcommand\wg{^{\wedge}}

\newcommand\ol{\overline}

\newcommand\ssim{/\!\!\sim}

\newcommand\End{\operatorname{End}}

\newcommand\Tr{\operatorname{Tr}\,}

\newcommand\qex{_{\operatorname{ex},q}}

\newcommand\id{\operatorname{id}}
\newcommand\lp{\operatorname{\!\langle\!}}
\newcommand\rp{\operatorname{\!\rangle\!}}

\newcommand{\isom}{\,\raise2pt\hbox{$\underrightarrow{\sim}$}\,}
\numberwithin{equation}{section}

\newtheorem{thm}{Theorem}[section]
\newtheorem{lem}[thm]{Lemma}
\newtheorem{cor}[thm]{Corollary}
\newtheorem{prop}[thm]{Proposition}

\def \para#1{\par\medskip\textbf{#1}
              \addtocounter{thm}{1}}
 
\begin{document}
\setlength{\baselineskip}{4.9mm}
\setlength{\abovedisplayskip}{4.5mm}
\setlength{\belowdisplayskip}{4.5mm}
%%%
%%%
\renewcommand{\theenumi}{\roman{enumi}}
\renewcommand{\labelenumi}{(\theenumi)}
\renewcommand{\thefootnote}{\fnsymbol{footnote}}
%%%
%\NoBlackBoxes
\parindent=20pt
\medskip
\title{Green functions associated to \\  complex reflection groups, II}
\author{Toshiaki Shoji}
\maketitle
\pagestyle{myheadings}
\markboth{T. SHOJI}{GREEN FUNCTIONS}
\begin{center}
Department of Mathematics  \\
Science University of Tokyo \\ 
Noda, Chiba 278-8510, Japan
\end{center}
\par\medskip
\begin{abstract}
Green functions associated to 
complex reflection groups $G(e,1,n)$ were discussed in the author's
previous paper.  In this
paper, we consider the case of complex reflection groups 
$W = G(e,p,n)$. Schur functions and Hall-Littlewood functions associated
to $W$ are introduced, and Green functions are described as the
transition matrix between those two symmetric functions.  
Furthermore, it is shown that these Green functions 
are determined by means of
Green functions associated to various $G(e',1,n')$.  Our result
 involves, as a special case, a combinatorial approach to the 
Green functions of type $D_n$.  
\end{abstract}
%%%%%%%%%%%%%%%%%%%%
%%%%%%%%%%%%%%%%%%%%
\bigskip
\medskip
\addtocounter{section}{-1}
\section{Introduction}
This paper is a continuation of [S].  In [S], Hall-Littlewood
functions associated to the complex reflection group $G(e,1,n)$
were introduced. Green functions associated to $G(e,1,n)$ are defined
as a solution of a certain matrix equation arising from the
combinatorics of $e$-symbols. It was shown that such Green functions
are obtained as coeffcients of the expansion of Schur functions in 
terms of Hall-Littlewood functions.  In the case where $e=2$, 
$G(e,1,n)$ coincides with the Weyl group of type $B_n$, and the Green
function in that case coincides with the Green function associated to
finite classical groups $Sp_{2n}(\Fq)$ or $SO_{2n+1}(\Fq)$ introduced 
by Deligne-Lusztig, in a geometric way. So our resut is
regarded as a first step towards the combinatorial description of such 
Green functions, just as in the case of Green polynomials of 
$GL_n(\Fq)$.
\par
In this paper, we take up the complex reflection group $G(e,p,n)$, and
show that a similar formalism as in the case of $G(e,1,n)$ works also 
for such groups. In the case of $G(e,p,n)$, symmetric functions such 
as Schur functions, Hall-Littlewood functions, etc.  appear as 
$p$-tuples of similar functions associated to the various 
complex reflection groups $G(e',1,n')$. In particular, Green functions 
associated to $G(e,p,n)$ can be described in terms of Green functions
associated to $G(e',1,n')$.  In the case where $e = p = 2$, the group
$G(e,p,n)$ is equal to the Weyl group of type $D_n$.  In this case
our Green functions coincide with the Green functions associated to
the finite groups $SO_{2n}(\Fq)$ of split type or non-split type.    
So our result implies, in this case, that the Green functions of type 
$D_n$ can be described completely in terms of various 
``Green functions'' of type $B_{n'}$.  However, note that the Green
function of type $B_{n'}$ appearing in this context is not the Green
function associated to $SO_{2n'+1}$.  They are the functions 
introduced in [S], associated to different type of symbols. 
\par
In the case where $n = 2$, the group $G(e,e,n)$ is equal to the
dihedral group of degree $2e$.  In this case, our Green function
coincides with the function obtained by Geck-Malle [GM] 
in connection with special pieces and unipotent characters of 
finite reductive groups (see [L2]).
\par
The author is grateful to Gunter Malle for some explicit computation 
by using a computer.  The tables in section 5 are due to him. 
\par\bigskip
\begin{center}
{\sc Table of contents}
\end{center}
\par\medskip
0. \ Introduction
\par
1. \ Schur-Weyl reciprocity for $G(e,p,n)$  
\par
2. \ Frobenius formula for $G(e,p,n)$  
\par
3. \ Hall-Littlewood functions associated to $G(e,p,n)$ 
\par
4. \ Green functions associated to $G(e,p,n)$ 
\par
5. \ Examples  
%%%%%
%%%%%
\par\bigskip
\section{Schur-Weyl reciprocity for $G(e,p,n)$}
\para{1.1.}
Let $\wt W \simeq \FS_n \ltimes (\ZZ/e\ZZ)^n$ be the imprimitive 
complex reflection group $G(e,1,n)$ acting on the complex vector space
$\CC^n$.  Let $e_1, \dots, e_n$ be the standard basis of $\CC^n$.
Then $\wt W$ is realized as a subgroup of $GL(\CC^n)$ consisting of $w$ 
such that $w(e_i) = \xi_ie_{x(i)}$ for $i = 1, \dots, n$, where 
$x \in \FS_n$ and $\xi_i$ is an $e$-th root of unity (depending on $w$).  
We fix a primitive $e$-th root of unity $\z$.  Then $\xi_i = \z^{a_i}$
with $a_i \in \ZZ/e\ZZ$, and $w \in \wt W$ can be written as
$w = x(a_1, \dots, a_n)$ uniquely.
The complex reflection group
$W_0 = G(e,e,n)$ is defined as the subgroup  of 
$\wt W$ consisting of $w \in \wt W$ such that 
$\sum_{i=1}^na_i \equiv 0 \pmod e$.
$W_0$ is 
a normal subgroup of $\wt W$ of index $e$.  Let us define, 
for $k = 0, \dots, e-1$,  $t_k \in \wt W$ by 
\begin{equation*}
t_k(e_i) = \begin{cases}
           \z e_k &\quad\text{ if } i = k, \\
             e_i &\quad\text{ otherwise.}
         \end{cases}
\end{equation*}  
Then $w = x(a_1, \dots, a_n)$ can be written as 
$w = x\prod_{i=1}^nt_i^{a_i}$.  Hence if we
put $\s = t_1$, $\wt W$ is a semidirect product of $W_0$ with the
cyclic group of order $e$ generated by $\s$.  For each factor $p$ of $e$, we
put $W = W_0\rtimes\lp \s^p\rp$. Then $W$ is the subgroup of
$\wt W$ consisting of $w$ such that 
$\sum_{i=1}^n a_i \equiv 0 \pmod p$, and is isomorphic to
the complex reflection group $G(e,p,n)$.
%%%%
%%%%
\para{1.2.}
For a factor $p$ of $e$, we put $d = e/p$.   
Let $\CP_{n,e}$ be the set of $e$-tuples 
$\Ba = (\a^{(0)}, \dots, \a^{(e-1)})$ of partitions such that 
$\sum_{k-0}^{e-1} |\a^{(k)}| = n$.  An element $\Ba \in \CP_{n,e}$ is
called an $e$-partition of $n$.  
We put $l(\Ba) = \sum_{k=0}^{e-1}l(\a^{(k)})$, where $l(\a^{(k)})$
denotes the number of parts of the partition $\a^{(k)}$.
Let us define an operator 
$\th = \th_p$ on $\CP_{n,e}$ by 
$\th(\Ba) = (\a^{(k-d)})$ for
$\Ba = (\a^{(k)}) \in \CP_{n,e}$.
We denote by $c_{\Ba,p}$ the number of elements in the $\th$-orbit
of $\Ba$.
Then the irreducible representations of $W$ are
described as follows. It is known that the set of isomorphism classe
of irreducible representations of $\wt W$ is parametrized by 
$\CP_{n,e}$.  Let $Z_{\Ba}$ be an irreducible $\wt W$-module 
corresponding to $\Ba \in \CP_{n,e}$.  Then $Z_{\Ba}$ is decomposed,
as a $W$-module, as
\begin{equation*}
\tag{1.2.1}
Z_{\Ba} = Z^p_{\Ba,1}\oplus\cdots\oplus Z^p_{\Ba,r},
\end{equation*} 
where $r = p/c_{\Ba,p}$, and $\s$ permutes factors in a cyclic way.
Each $Z^p_{\Ba,i}$ is an irreducible $W$-module, mutually 
non-isomorphic of the same dimension.   
Furthermore, if 
$\Ba$ and $\Bb$ are in the same orbit under $\th$, then 
$Z_{\Ba} \simeq Z_{\Bb}$ as $W$-modules.  We denote by 
$\x^{\Ba}$ the irreducible character of $\wt W$ afforded by $Z_{\Ba}$,
and by $\x_p^{\Ba,i}$
the irreducible character of $W$ afforded by $Z^p_{\Ba,i}$ for 
$i = 1, \dots, p/c_{\Ba,p}$. 
We write $\Ba \sim_p \Bb$ if $\Ba$ and $\Bb$ are in the same
$\th$-orbit, and by $\CP_{n,e}\ssim_p$ the set of 
$\th$-orbits  in $\CP_{n,e}$.  Then the set 
$W\wg$  of irreducible characters  of $W$ is 
parametrized as
\begin{equation*}
\tag{1.2.2}
W\wg =  \{ \x_p^{\Ba,i} \mid \Ba \in 
               \CP_{n,e}\ssim_p, 1 \le i \le p/c_{\Ba,p} \}.
\end{equation*} 
\para{1.3.}
Here we recall the Schur-Weyl reciprocity between $\wt W$ and 
a cetaiin Levi subgroup of a general linear group over $\CC$.
Let $V = \bigoplus_{i=0}^{e-1} V_i$ be a vector space over $\CC$
with $\dim V_i = m_i$.  We fix a basis 
$\CE = \{ v_j^{(k)} \mid 1 \le j \le m_k\}$ of $V_k$ for $0 \le k \le e-1$.
Then 
$v_1^{(0)}, \dots, v_{m_0}^{(0)}, v_1^{(1)}, \dots, v_{m_1}^{(1)},
\dots$
gives a basis of $V$, which we write in this order as 
$v_1, \dots, v_M$ with $M = \sum m_i$.
Let $G = GL_{m_0} \times \cdots \times GL_{m_{e-1}}$.
Here $GL_{m_i}$ acts on $V_i$ in a natural way.  Hence we have an action
of $G$ on $V$, and so on the $n$-fold tensor space $V^{\otimes n}$.  
On the other hand, $\FS_n$ acts on $V^{\otimes n}$ by permuting the 
factors of the tensor product.  This action commutes with the action
of $G$.  We extend the action of $\FS_n$ to that of $\wt W$ as follows:
For a basis element $v = v_{i_1}\otimes\cdots\otimes v_{i_n}$ in $V$, 
we put 
\begin{equation*}
t_k(v) = \z^jv 
\end{equation*} 
if $v_{i_k} \in V_j$.  Then this action of $t_1, \dots, t_n$ on 
$V$ gives rise to an action of $\wt W$, commuting with the action of $G$. 
It is known that the following Schur-Weyl reciprocity holds (e.g., see
[SS], where the Hecke algebra version is discussed).
\par\medskip\noindent
(1.3.1)  Let $\s_1: \CC\wt W \to \End V^{\otimes n}, 
\r_1: \CC G \to \End V^{\otimes n}$ be the representations of 
$\wt W$ and $G$, respectively, ($\CC\wt W$, etc.  denote the group algebra
of $\wt W$, etc.  over $\CC$).  Then $\s_1(\CC\wt W)$ and 
$\r_1(\CC G)$ are the centralizer algebras of each other in 
$\End V^{\otimes n}$. More precisely, the following holds.  
Put $\Bm = (m_0, \dots, m_{e-1})$, and let $\vL_{\Bm}$ 
be the set of 
$\Ba \in \CP_{n,e}$ such that $l(\a^{(i)}) \le m_i$.  
Then  
$\wt W \times G$-module $V^{\otimes n}$ is decomposed as
\begin{equation*}
V^{\otimes n} = 
  \bigoplus_{\Ba \in \vL_{\Bm}}
      Z_{\Ba}\otimes V_{\Ba},
\end{equation*}
where $Z_{\Ba}$ is the irreducible $\wt W$-module as before, and
$V_{\Ba}$ is an irreducible $G$-module with highest weight $\Ba$.
In particular, if $m_i \ge n$ for any $i$,  all the irreducible 
$\wt W$-modules are realized in $V^{\otimes n}$. 
%%%%
%%%%
\para{1.4.}
We shall extend the Schur-Weyl reciprocity to the case of $W$.
Here we pose the following assumption on $\Bm$.
\par\medskip\noindent
(1.4.1) \  $m_k = m_{k+d}$  for $k = 0, 1, \dots, p-1$ in
$\Bm =(m_0, \dots, m_{e-1})$.
\par\medskip  
Let us define a linear automorphism $\t = \t_p$ on $V$ by
$\t(v_j^{(k)}) = v_j^{(k+d)}$ for $1 \le j \le m_k$ and
 for integers $k$
(here we regard $k \in \ZZ/e\ZZ$).  Then $\t$ is an element of $GL_M$
of order $p$, and normalizes the subgroup $G$.  We denote by $\wt G$
the subgroup of $GL_M$ generated by $G$ and $\t$, which is isomorphic
to the semidirect product $G \rtimes \lp \t\rp$, where $\t$ acts on 
$G = GL_{m_0} \times\cdots\times  GL_{m_{e-1}}$ by 
$\t(GL_{m_k}) = GL_{m_{k+d}}$.   Now we have an action of 
$\wt G$ on $V^{\otimes n}$, where
 $\t$ acts on $V^{\otimes n}$ by
\begin{equation*}\
\t(v_{i_1}^{(k_1)}\otimes\cdots\otimes v_{i_n}^{(k_n)})
  = v_{i_1}^{(k_1+d)}\otimes\cdots\otimes v_{i_n}^{(k_n+d)}.
\end{equation*} 
We denote this action also by $\r_1$.
Here $\r_1(\t)$ is an automorphism of order $p$ on $V^{\otimes n}$, 
commuting with the action of $\FS_n$.
Moreover, we have $\t \s\t\iv = \z^{-d} \s$ on $V^{\otimes n}$, and so 
$\r_1(\t)$ normalizes the subalgebra $\r_1(\CC \wt W)$ of 
$\End V^{\otimes n}$.  
It is easy to check, for $w \in \wt W$,  that $\s_1(w)$ commutes 
with $\r_1(\t)$ if and only if
$w \in W$.  Thus we have an action of 
$W \times \wt G$ on $V^{\otimes  n}$.
We consider the decomposition of $V^{\otimes n}$ as in (1.3.1).
Since $\t(V_{\Ba}) \isom V_{\th(\Ba)}$, $\t$ maps 
$Z_{\Ba} \otimes V_{\Ba}$ onto $Z_{\th(\Ba)}\otimes V_{\th(\Ba)}$. 
Here we may choose, as a model of $Z_{\Ba}$ in $V^{\otimes n}$,
the space of highest weight vectors with highest weight $\Ba$ in 
$V^{\otimes n}$.  Then we have $\t(Z_{\Ba}) = Z_{\th(\Ba)}$.
Put $\vG = <\t> \simeq \ZZ/p\ZZ$, and let $\vG_{\Ba}$ be 
the stabilizer of $Z_{\Ba}$ in $\vG$. 
Then $Z_{\Ba} = \lp \t^c\rp$ with $c = c_{\Ba,p}$, and 
$Z_{\Ba}$ turns out to be a $\vG_{\Ba}$-module.
For each $\f \in \vG_{\Ba}\wg$, we denote by
$Z_{\Ba,\f}$ the $\f$-isotropic subspace of $Z_{\Ba}$, which
is stable by $W$.
Now $\s$ permutes factors $Z_{\Ba,\f}$ transitively, 
and $|\vG_{\Ba}| = p/c_{\Ba,p}$.  It follows that 
each $Z_{\Ba,\f}$ becomes an irreducible $W$-module and
\begin{equation*}
Z_{\Ba} = \bigoplus_{\f \in \vG_{\Ba}\wg} Z_{\Ba,\f}
\end{equation*}
gives the decomposition of $Z_{\Ba}$ into irreducible
$W$-modules given in (1.2.1).
On the other hand, let $V'_{\Ba,\f}$ be a $G$-submodule of 
$V^{\otimes n}$ generated by a highest weight vector in $Z_{\Ba,\f}$.
Then $V'_{\Ba,\f}$ is isomorphic to $V_{\Ba}$, and is stable under
the action of $\t^c$.  
We denote by $V_{\Ba,\f}$ the $\wt G$-submodule of $V^{\otimes n}$ 
generated by $V'_{\Ba,\f}$.
Then $V_{\Ba,\f}$ is isomorphic to the $\wt G$-module 
induced from $\vG_{\Ba}\ltimes G$-module $V'_{\Ba,\f}$, and 
turns out to be an irreducible $\wt G$-module.  Moreover, 
$V_{\Ba,\f}$ are mutually non-isomorphic for distinct pairs 
$(\Ba,\f)$.
Put 
\begin{equation*}
\vL^p_{\Bm} = 
 \{ (\Ba,\f) \mid \Ba \in \vL_{\Bm}\ssim_p, 
                     \f \in \vG_{\Ba}\wg \},
\end{equation*}
where $\vL_{\Bm}\ssim_p$ denotes the set of $\th$-orbits in $\vL_{\Bm}$.
It follows from the above discussion, we have the following 
Schur-Weyl reciprocity between $W$ and $\wt G$.
\begin{prop} %%%Poro. 1.5
$\s_1(\CC W)$ and $\r_1(\CC\wt G)$ are mutually the full centralizer
algebras of each other in $\End V^{\otimes n}$.  Moreover, 
$W \times \wt G$-module $V^{\otimes n}$ is decomposed as
\begin{equation*}
V^{\otimes n} = \bigoplus_{(\Ba,\f)  \in \vL^p_{\Bm}}
                Z_{\Ba,\f}\otimes V_{\Ba,\f}.
\end{equation*}
\end{prop}
%%%%
\section{Frobenius formula for $G(e,p,n)$}
%%%%%
\para{2.1.}
In the remainder of this paper, we assume that 
$\Bm = (m_0, \dots, m_{e-1})$ satisfies the condition 
(1.4.1) and the condition
that $m_i \ge n$ for $i = 0, \dots, e-1$. 
Then 
any irreducicible $W$-module is realized as
$Z_{\Ba, \f}$ by Proposition 1.5. In this case, $\vL_{\Bm}^p$
coincides with the set 
\begin{equation*}
\wt\CP_W = \{ (\Ba, \f) \mid \Ba \in \CP_{n,e}\ssim_p, 
                   \f \in \vG_{\Ba}\wg\}.
\end{equation*}
We denote by $\x^{\Ba,\f}$ the irreducible character of $W$
corresponding to the $W$-module $Z_{\Ba,\f}$.  Under this
notation, the parametrization of $W\wg$ in (1.2.2) can be 
modified as  
\begin{equation*}
\tag{2.1.1}
W\wg = \{ \x^{\Ba,\f} \mid (\Ba, \f) \in \wt\CP_W\}.
\end{equation*}
\par 
For later use, we consider a more general situation.  Assume that
$q$ is a factor of $e$ and that 
$\lp \s^{q}\rp \times \lp \s^p \rp$ is a subgroup of 
$\lp \s\rp \simeq \ZZ/e\ZZ$, i.e., $e/q$ and $e/p$ are 
prime each other.  The typical case is that $W = W_0$ and 
$q$ is any factor of $e$.  Now one can form a semidirect product
$\lp \s^q\rp\ltimes W$ as a subgroup of $\wt W$. 
We are now 
interested in the set of $\s^q$-stable characters in  $W\wg$. 
It follows from the discussion 
in 1.4 that the irreducible 
$W$-module $Z_{\Ba,\f}$ is $\s^{q}$-stable if and only if $q$ is 
divisible by $|\vG_{\Ba}| = p/c_{\Ba,p}$.  So we put
\begin{equation*}
\wt\CP_W^q = \{ (\Ba,\f) \in \wt\CP_W \mid 
    qc_{\Ba,p}\equiv 0 \!\!\pmod p \}.
\end{equation*}
We also write it as $\wt\CP_{n,e,p}^q$ to make the dependence 
on $n,e,p$ more explicit.
Let  
$W\wg\qex$ be the set of $\s^q$-stable irreducible 
characters of $W$. Thus $W\wg\qex$ is given as 
\begin{equation*}
\tag{2.1.2}
W\wg\qex =  \{ \x^{\Ba,\f} \in W\wg \mid (\Ba,\f) \in \wt\CP^q_W \}.
\end{equation*} 
%%%%
\para{2.2.}  
As in [S], we identify $\CP_{n,e}$ with the set 
$Z_n^{0,0} = Z_n^{0,0}(\Bm)$.   
Here $Z_n^{0,0}$ is the set of $e$-partitions $\Ba$ written as 
$\Ba = (\a^{(0)}, \dots, \a^{(e-1)})$ with a partition  
$\a^{(k)}: \a^{(k)}_1 \ge \cdots \ge \a^{(k)}_{m_k} \ge 0$. 
We prepare the indeterminates 
$x_i^{(k)}$ for $0 \le k <e, 1 \le i \le m_k$. We write 
$x = \{x_i^{(k)}\}$, and also write as 
$x^{(k)} = \{x_1^{(k)}, \dots, x_{m_k}^{(k)}\}$ for a fixed $k$.
Recall that a power sum symmetric function $p_{\Ba}(x)$ is defined 
for $\Ba \in Z_n^{0,0}$ as  
\begin{equation*}
\tag{2.2.1}
p_{\Ba}(x) = \prod_{k=0}^{e-1}\prod_{j=1}^{m_k}p^{(k)}_{\a^{(k)}_j}(x).
\end{equation*}
Here we put, for each integer $r \ge 1$,  
\begin{equation*}
\tag{2.2.2}
p_r^{(i)}(x) = \sum_{j=0}^{e-1}\z^{ij}p_r(x^{(j)}),
\end{equation*}
with usual $r$-th powersum symmetric functions $p_r(x^{(j)})$ with
variables $x^{(j)}$, and put $p_0^{(i)}(x) = 1$.
Also, Schur functions $s_{\Ba}(x)$ and monomial symmetric functions
$m_{\Ba}(x)$ are defined as
\begin{equation*}
s_{\Ba}(x) = \prod_{k=0}^{e-1}s_{\a^{(k)}}(x^{(k)}), \qquad
m_{\Ba}(x) = \prod_{k=0}^{e-1}m_{\a^{(k)}}(x^{(k)}),
\end{equation*}
by using usual Schur functions $s_{\a^{(k)}}$ and 
monomial symmetric functions $m_{\a^{(k)}}$ associated to 
partitions $\a^{(k)}$.
\par
It is known that the set of conjugacy classes in $\wt W$ is in
bijection with $\CP_{n,e}$. The explicit correspondence will be
given later in 2.3.  We denote by $w_{\Bb}$ a representative of 
the conjugacy
class in $\wt W$ corresponding to $\Bb \in \CP_{n,e}$.
Then the Frobenius formula for the irreducible characters 
of $\wt W$ is given as follows.
\par\medskip\noindent
(2.2.3) \ (Frobenius formula for $\wt W$).  Let $\Bb \in \CP_{n,e}$.
Then we have
\begin{equation*}
p_{\Bb} = \sum_{\Ba \in \CP_{n,e}}\x^{\Ba}(w_{\Bb})s_{\Ba}.
\end{equation*} 
%%%%
\para{2.3.}
We want to generalize (2.2.3) to the case of irreducible characters on
$\s^qW$.  
First of all we describe the conjugacy classes in $\wt W$ and 
$W$ more precisely.
%%%%
The correspondence between the conjugacy classes in $\wt W$ and 
$\CP_{n,e}$ is given explicitly as follows.
Assume that $w$ maps $e_{i_1}, e_{i_2}, \dots, e_{i_m}$ in a cyclic
way, up to scalar, and leaves other $e_j$ unchanged.  
If $w^m(e_{i_1}) = \z^ke_{i_1}$, 
we say that $w$ is a $k$-cycle $(i_1, \dots, i_m)$ of length $m$. 
Now $w \in \wt W$ can be written as a product of various disjoint 
cycles.
If we group up, for a fixed $k$, the $k$-cycles among them, 
it produces a partition $\a^{(k)}$, 
and $\Ba = (\a^{(0)}, \dots, \a^{(e-1)})$ turns out to be
an $e$-partition of $n$.  This gives the required bijection.  
We denote by $C_{\Ba}$ the conjugacy class of $\wt W$ corresponding to
$\Ba \in \CP_{n,e}$.
\par
For $\Ba = (\a^{(0)}, \dots, \a^{(e-1)}) \in \CP_{n,e}$, we define
an integer $\vD(\Ba)$ by 
\begin{equation*}
\tag{2.3.1}
\vD(\Ba) = \sum_{k = 0}^{e-1}l(\a^{(k)})k.
\end{equation*}
It is easily checked taht the conjugacy class
$C_{\Ba}$ belongs to $W$ if and only if $\vD(\Ba) \equiv 0 \pmod p$.
The class $C_{\Ba}$ is decomposed into several conjugacy classes in 
$W$.   
\par
More generally, we consider a coset $\s^qW$ in $\wt W$.
Then $\s^qW$ is invariant under the adjoint action of $\wt W$, 
and the class $C_{\Ba}$ lies in $\s^qW$ if and only if 
$\vD(\Ba) \equiv q \pmod p$.
The conjugacy class $C_{\Ba} \subset \s^qW$ is also decomposed 
into several $W$-orbits. The complete description of such 
$W$-orbits will be
given in Proposition 2.13.  
\para{2.4.}
Let $D$ be an operator on $V^{\otimes n}$ defined by
\begin{equation*}
D(v_{i_1}\otimes v_{i_2}\otimes\cdots\otimes v_{i_n}) 
     = x_{i_1}x_{i_2}\cdots x_{i_n}
          (v_{i_1}\otimes v_{i_2}\otimes\cdots\otimes v_{i_n}).
\end{equation*}
(Although we must replace $V$ by the scalar extension $K\otimes_{\CC} V$ 
with $K = \CC(x_j^{(k)})$, we use the same notation 
by abbreviation).  Then $D$ commutes with the action of $\wt W$, and
the Frobenius formula (2.2.3) is obtained by computing the trace
$\Tr(Dw, V^{\otimes n})$ for $w \in \wt W$ in two different ways. 
We follow this strategy to establish the Frobenius formula for 
the characters of $W$, in our case by computing the traces
$\Tr(D\t^jw, V^{\otimes n})$ for $j = 0, \dots, p-1$.
\par
For an integer $j \, (0 \le j < p)$, let $h_j$ be order 
of the image $\bar j$ of $j$ in $\ZZ/p\ZZ$. Then $p$ is written as 
$p = h_jj_1$,  where
$j_1$ is the greatest common divisor of $j$ and $p$.   
We introduce  new indeterminates $\CX_j = \{X_i^{(k)}\}$ by
\begin{equation*}
\tag{2.4.1}
X_i^{(k)} = x_i^{(k)}x_i^{(k+jd)}x_i^{(k+2jd)}\cdots x_i^{(k+ (h -1)jd)}
   \qquad ( 0 \le k < j_1d, 1 \le i \le m_k), 
\end{equation*} 
with $h = h_j$.
We write also as 
$\CX_j^{(k)} = \{ X_1^{(k)}, \dots, X_{m_k}^{(k)}\}$ for a fixed $k$.
One can define a  function 
$p_r^{(i)}(\CX_j)$  
with respect to the variables 
$\CX_j = \{ \CX_j^{(0)}, \dots, \CX_j^{(j_1d-1)}\}$ by 
modifying (2.2.2) as follows.  
\begin{equation*}
\tag{2.4.2}
p_r^{(i)}(\CX_j) = \sum_{k=0}^{j_1d-1}\z^{ik}p_r(\CX_j^{(k)}).
\end{equation*} 
\par
Then we have the following lemma.
%%%%
\begin{lem}%%%% Lemma 2.5
Let $w \in \wt W$ be an $f$-cycle of length $n$.  Assume that  
$w$ is of the form  
$w = t_1^at_n^bz$, where $z$ is a cyclic permutation 
$(1,2,\dots,n) \in \FS_n$ and $a+b \equiv f \pmod p$. 
Then, for each $j$ such that $0 \le j < p$,
\begin{equation*}
\Tr(D\t^jw, V^{\otimes n}) = 
                \begin{cases}
          h\z^{-(f + b)jd}p_{n/h}^{(f)}(\CX_j)
        &\quad\text{ if  $h\mid n$ and $h \mid f$,} \\
                     0 
        &\quad\text{ otherwise,}
                \end{cases}
\end{equation*}
where $h = h_j$ as before.
\end{lem}
\begin{proof}
We have
\begin{equation*}
\t^jw(v_{i_1}^{(k_1)}\otimes v_{i_2}^{(k_2)}
      \otimes\cdots\otimes v_{i_n}^{(k_n)})
          = \z^{ak_n + bk_{n-1}}v_{i_n}^{(k_n + jd)}\otimes
               v_{i_1}^{(k_1 + jd)}\otimes
             \cdots\otimes v_{i_{n-1}}^{(k_{n-1} + jd)}.
\end{equation*}
Now assume that $\Tr(D\t^jw, V^{\otimes n}) \ne 0$.  Then there exists 
$v^{(k_1)}_{i_1}, \dots, v_{i_n}^{(k_n)} \in \CE$ satisfying 
the relation 
\begin{equation*}
v_{i_s}^{(k_s)} = v_{i_{s-1}}^{(k_{s-1} + jd)}\qquad 
\text{ for } s = 1, \dots, n. 
\end{equation*}
(Here we regard $v_{i_0} = v_{i_n}$, and $k_0 = k_n$).  This
implies that 
\begin{align*} 
\tag{2.5.1}
v_{i_1}^{(k_1)} &= v_{i_1}^{(k_1 + njd)},  \\ 
v_{i_s}^{(k_s)} &= v_{i_1}^{(k_1 + (s-1)jd)} \quad\text{ for }
s = 2, \dots, n.
\end{align*}
In particular, we have $njd \equiv 0 \pmod e$, or 
equivalently $nj \equiv 0 \pmod p$.
\par
Under the notation in (2.4.1), we have 
\begin{equation*}
\tag{2.5.2}
\Tr(D\t^jw, V^{\otimes n}) = \z^{-(a+2b)jd}
   \sum_{k=0}^{j_1d-1}\sum_{s=0}^{h-1}\z^{f(k+sjd)}
          \sum_{i=1}^{m_k}(X_i^{(k)})^{n/h}.
\end{equation*}
This implies that $\Tr(D\t^jw, V^{\otimes n}) = 0$ unless
$fj \equiv 0 \pmod p$.  Here note that the condition 
$nj \equiv 0 \pmod p$ is equivalent to $h \mid n$, and similarly
for $fj \equiv 0 \pmod p$.
The second assertion follows from this.  
\par
We now assume that $nj \equiv 0$ and $fj \equiv 0 \pmod p$.
Then by (2.4.2) and (2.5.2), we have
\begin{equation*}
\Tr(D\t^jw, V^{\otimes n}) = h\z^{-(a+2b)jd}p_{n/h}^{(f)}(\CX_j).
\end{equation*}
This shows the first assertion, and the lemma is proved.
\end{proof}
\para{2.6.}
We consider the general case where $w$ is in the class
$C_{\Ba}$ in $\wt W$.
We fix $j$ and let $h = h_j$ and $j_1$ be as before. 
For $\Ba = (\a^{(k)}_i)$, we consider the following 
condition (for a fixed $j$). 
\begin{align*}
\tag{2.6.1}
&h \mid \a^{(k)}_i  \quad\text{ for any }\a_i^{(k)},
                \text{ and} \\
&h \mid k \quad\text{ if } |\a^{(k)}| \ne 0. 
\end{align*} 
Assume that $\Ba \in \CP_{n,e}$ 
satisfies (2.6.1).  We define a $j_1d$-partition 
$\Bb = (\b^{(0)}, \dots, \b^{(j_1d-1)})$ of $n/h$ by
$\b^{(k/h)}_i = \a^{(k)}_i/h$ for any $\a^{(k)}_i$ such that
$h \mid k$.  
We write $\Bb$ as $\Bb = \Ba[j]$.  Then we have 
$h\vD(\Ba[j]) = \vD(\Ba)$.  
Note that in this situation (2.4.2) can be regarded as a formula
analogous to (2.2.2) for variables $\CX_j$, by replacing $\z$ by
$\z^h$.  Thus one can define $p_{\Ba[j]}(\CX_j)$ just as in (2.2.1)
as a product of various $p_{\b_i^{(k)}}(\CX_j)$. 
\par
As in 2.3, $w$ can be written as a product of various cycles
associated to $\Ba$. Under the conjugation in $W$, $w$ is
changed to an element $w'$ of the following type; the cycle
corresponding to $\a_i^{(k)}$ is of the form 
$t_{i_1}^kz$, where $z = (i_1, i_2, \dots, i_m) \in \FS_m$ with 
$m = \a_i^{(k)}$, except one cycle.  The excepted one is of the form 
$t_{i_1}^at_{i_m}^bz$, where 
$z = (i_1, i_2, \dots, i_m)$ with $a + b \equiv k\pmod e$. 
We may further assume that each cyclic permutation $z$ occuring in
the decomposition of $w'$ is of the form 
$z = (i_1, i_1+1, \dots, i_1+(m-1))$.  We write this element as 
$w' = w_{\Ba}(b)$.
\par
Note that $\Tr(D\t^jw, V^{\otimes n})$ remain unchanged if $w$ is
replaced by its conjugate $w'$, since $W$ commutes with 
$D$ and $\t$.   Then $\Tr(D\t^jw', V^{\otimes n})$ can be computed 
by making use of Lemma 2.5.  So, as a corollary to Lemma 2.5, we have
%%%%
\begin{prop}%%%% Prop.2.7.
Assume that $w \in C_{\Ba}$ in $\wt W$, 
and that $w$ is conjugate under $W$ to $w_{\Ba}(b)$.
Then $\Tr(D\t^jw, V^{\otimes n}) = 0$ unless $\Ba$ satisfies the
condition (2.6.1).
If $\Ba$ satisfies the condition (2.6.1), then
\begin{equation*}
\Tr(D\t^jw, V^{\otimes n}) = 
     h^{l(\Ba[j])}\z^{-(\vD(\Ba)+b)jd}p_{\Ba[j]}(\CX_j).
\end{equation*}
\end{prop}  
%%%%%
\para{2.8.}
Next, we consider the decomposition    
$V^{\otimes n} = \bigoplus V_{\Ba,\f}\otimes Z_{\Ba,\f}$ as 
given in Proposition 1.5.  
Assume that $w \in \s^qW$.  
If $\x^{\Ba,\f} \in W\wg$ is not $\s^q$-stable,
then $D\t^jw$ maps $V_{\Ba,\f}\otimes Z_{\Ba,\f}$ onto a different
factor.  It follows that
\begin{equation*}
\tag{2.8.1}
\Tr(D\t^jw, V^{\otimes n}) = 
     \sum_{(\Ba,\f) \in \wt\CP^q_W}
          \Tr(D\t^jw, V_{\Ba,\f}\otimes Z_{\Ba,\f}).
\end{equation*}
\par
Let $\x^{\Ba,\f}$ be a $\s^q$-stable irreducible
character of $W$. Then $\x^{\Ba,\f}$ can be extended to an
irreducible character of $\lp \s^q \rp \ltimes W \simeq W_{e,p',n}$, 
where $p' = pq/e$.  Let $\vG' \simeq \ZZ/p'\ZZ$ 
be the subgroup of $\vG \simeq \ZZ/p\ZZ$ generated by 
$\t^{e/q}$, 
and $\vG'_{\Ba}$ the stabilizer of $Z_{\Ba}$ in $\vG'$.  
Since $q$ is divisible by $|\vG_{\Ba}|$ by (2.1.2) and
since $(e/p,e/q)=1$, we see 
that $\vG'_{\Ba}$ coincides with $\vG_{\Ba}$.  Hence Irreducible
characters of $W$ occuring in the decomposition of $\x^{\Ba}$
are also regarded as irreducible characters of $W_{e,p',n}$.  
In this way, we can fix an extension $\wt\x^{\Ba,\f}$
of $\x^{\Ba,\f}$ to $W_{e,p',n}$. 
Accordingly, the irreducible $W$-module  $Z_{\Ba,\f}$ is
extended to the irreducible $W_{e,p',n}$-module corresponding to 
$\wt\x^{\Ba,\f}$, which we denote by $\wt Z_{\Ba,\f}$. 
\par
Remember that $V_{\Ba,\f}\otimes Z_{\Ba,\f}
     \simeq \bigoplus_{\Bb \in O(\Ba)}V'_{\Bb,\f}\otimes Z_{\Bb,\f}$
as $G \times W$-module, and $\t$ permutes the summands of the right
hand side. In the case where $c_{\Ba,p} \mid j$, $\t^j$ leaves each
summand $V_{\Bb,\f}'\otimes Z_{\Bb,\f}$ stable.
The last  space is extended to a $W_{e,p',n}$-module 
$V'_{\Bb,\f}\otimes\wt Z_{\Bb,\f}$, and we have an action of 
$w$ on 
$V_{\Ba,\f}\otimes\wt Z_{\Ba,\f} 
  \simeq \bigoplus_{\Bb \in O(\Ba)} V'_{\Bb,\f}\otimes\wt Z_{\Bb,\f}$.  
\par
Assume that $\Ba \in \CP_{n,e}$ satisfies the condition that 
$c_{\Ba,p} \mid j$.  We
put $\Ba\{j\} = (\a^{(0)}, \dots, \a^{(j_1d-1)})$. Since   
$j_1d$ is divisible by $c_{\Ba,p}$, $\Ba\{j\}$ turns out to be
a $j_1d$-partition of $n/h$. Also note that, $\t^j \in \vG_{\Ba}$.
Under this setting, We have the following lemma.
%%%%
\begin{lem} %%% Lemma 2.9
Let $\Ba \in \CP_{n,e}$, and $q$ be as in 2.1.  Assume that
$\x^{\Ba,\f}$ is $\s^q$-stable.  Then 
for $w \in \s^qW$, we have
\begin{align*}
\Tr(D\t^jw, &V_{\Ba,\f}\otimes\wt Z_{\Ba,\f}) = \\
       &\begin{cases}
     \f(\t^j)\displaystyle\sum_{0 \le i <c }\z^{qid}s_{\th^i(\Ba)\{j\}}
                   (\CX_j)\wt\x^{\Ba,\f}(w) 
         &\quad\text{ if $c \mid j$},  \\
              0 &\quad\text{ otherwise},         
       \end{cases}                        
\end{align*}
where $c = c_{\Ba,p}$, and  $s_{\Ba\{j\}}(\CX_j)$ etc. are 
Schur functions  given in 2.2 with
resepect to the variables $\CX_j$.
\end{lem}
%%%%
\begin{proof}
By the previous remark, $D\t^jw$ permutes the factors 
$V'_{\Bb,\f}\otimes Z_{\Bb,\f}$ if $j$ is not a multiple of
$c$, and so $\Tr(D\t^jw, V_{\Ba,\f}\otimes \wt Z_{\Ba,\f}) = 0$.
Now assume that $c\mid j$.  Then we can write
\begin{equation*}
\Tr(D\t^jw, V_{\Ba,\f}\otimes\wt Z_{\Ba,\f}) = 
         \sum_{\Bb \in O(\Ba)}\Tr(D\t^j, V'_{\Bb,\f})\wt\x^{\Bb,\f}(w).
\end{equation*}
Since $V'_{\Bb,\f} = \f\otimes V'_{\Bb,1}$ as 
$\lp\t^c\rp\ltimes G$-modules, we have
$\Tr(D\t^j, V'_{\Bb,\f}) = \f(\t^j)\Tr(D\t^j, V'_{\Bb,1})$.
We note that 
\begin{equation*}
\tag{2.9.1}
\Tr(D\t^j, V'_{\Bb,1}) = s_{\Bb\{j\}}(\CX_j). 
\end{equation*} 
\par
In fact, $V'_{\Bb,1}$ is an irreducible $G$-module 
isomorphic to $V_{\Bb}$, 
and it is known from (2.2.3) that $\Tr(D, V_{\Bb}) = s_{\Bb}(x)$.
$V_{\Bb}$ has a basis consisting of weight vectors, on which $D$ acts
diagonally.  Here for a vector 
$v = v_{i_1}^{(k_1)}\otimes\cdots\otimes v_{i_n}^{(k_n)} 
       \in V^{\otimes n}$, 
the weight of $v$ is given by an elemenet 
$\Bg(v) = (\g^{(0)}, \dots, \g^{(e-1)}) \in Z_n^{0,0}$, 
where $\g_i^{(k)}$ is the number of $v_i^{(k)}$ occuring in 
the expression of $v$. 
On the other hand, $\t^j$ acts on $V_{\Bb}$, by permuting the weight 
vectors; if $v$ is a weight vector with weight $\Bg$, then $\t^j(v)$
is also a weight vector with weight 
$\Bg' = (\g^{(-jd)}, \g^{(1-jd)}, \dots, \g^{(e-1 - jd)})$.
Hence, in order to compute $\Tr(D\t^j, V_{\Bb})$, we have only to 
consider the weight vectors in $V_{\Bb}$ whose weights are of the type
$\Bg = (\g_i^{(k)})$ such that $\g^{(k)} = \g^{(k+jd)}$ for each $k$. 
Note that monomials $\prod (x_i^{(k)})^{\g_i^{(k)}}$ obtained as
weights of $D$ produce the Schur function $s_{\Bb}(x)$.  Then the
corresponding weight for $D\t^j$ is given by  
$\prod (X_i^{(k)})^{\g_i^{(k)}}$ for $\Bg$ as above 
(the product is taken for $k$ such
that $0 \le k < j_1d$). 
Hence $\Tr(D\t^j, V'_{\b,1})$ is obtained by picking up the 
monomials of this type from 
$s_{\Bb}(x) = \prod_{k=0}^{e-1}s_{\Bb^{(k)}}(x^{(k)})$, which
coinicdes with $s_{\Bb\{ j\}}(\CX_j)$. 
This proves (2.9.1).
\par
To prove the lemma, it is enough to show that 
\begin{equation*}
\tag{2.9.2}
\wt\x^{\th^i(\Ba),\f}(w) = \z^{qid}\wt\x^{\Ba,\f}(w)
\end{equation*}
for $i = 0, \dots, c -1$.
We show (2.9.2).  Since $\wt\x^{\Ba,\f}$ and $\wt\x^{\th^i(\Ba),\f}$ are 
both extensions of $\x^{\Ba,\f}$, there exists a linear character 
$\vf$ of the cyclic group $\lp \s^q\rp$ such that 
$\wt\x^{\th^i(\Ba),\f} = \vf\otimes\wt\x^{\Ba,\f}$.
On the other hand, since $\t \s\t\iv = \z^{-d}\s$ on $V^{\otimes n}$, 
we see that the action of $\s^q$ on 
$\wt Z_{\th^i(\Ba),\f} = \t^i(\wt Z_{\Ba, \f})$ corresponds, under the
isomorphism $\t^i$ of $W$-modules,  to the
action of $\s^q$ on $\wt Z_{\Ba,\f}$ multiplied by $\z^{qid}$.
This implies (2.9.2), and the lemma follows.
\end{proof}  
\para{2.10.}
In order to formulate the Frobenius formula for the $\s^q$-stable 
characters of $W$,
we shall define Schur functions and power sum symmetric functions 
associated to $\s^qW$.
Let $(\Ba,\f) \in \wt\CP^q_W$, and put, for each
$0 \le j < p$,  
\begin{equation*}
\tag{2.10.1}
s_{\Ba,\f}^j(x) = \begin{cases}
                    \f(\t^j)\sum_{0 \le i < c}
            \z^{qid}s_{\th^i(\Ba)\{j\}}(\CX_j)
                  &\quad\text{ if } c \mid j, \\
                    0 &\quad\text{ otherwise},
                   \end{cases} 
\end{equation*} 
where $c = c_{\Ba,p}$.
We put $\Bs_{\Ba,\f}(x) = (s_{\Ba,\f}^j(x))_{0 \le j < p}$,
and call it Schur function for $\s^qW$ associated to 
$(\Ba,\f)$.
\par
Next, we consider power sum symmetric functions.  For a pair 
($\Ba,b)$, where $\Ba \in \CP_{n,e}$ and $0 \le b < e$, 
and for $0 \le j < p$, we put
\begin{equation*}
\tag{2.10.2}
p_{\Ba,b}^j(x) = \begin{cases} 
         h^{l(\Ba[j])}\z^{-(\vD(\Ba) +b)jd}p_{\Ba[j]}(\CX_j)
              &\quad\text{ if $\Ba$ satisfies (2.6.1)}, \\
          0 &\quad\text{ otherwise}.
                 \end{cases}
\end{equation*}
(For the notation, see 2.6.)  We put 
$\Bp_{\Ba,b}(x) = (p^j_{\Ba,b}(x))_{0 \le j < p}$, 
and call it a power sum symmetric
function for $\s^qW$ associated to the pair $(\Ba,b)$.
\par
We are now ready to formulate a Frobenius formula for the characters
on $\s^qW$ as follows. The proof is immediate from Proposition 2.7, 
(2.8.1) and Lemma 2.9.  
\begin{prop}[{\bf Frobenius formula for $\s^qW$}]%%%%Prop. 2.11
Let $(\Bb, b)$ be such that 
$\Bb \in \CP_{n,e}$ and $0 \le b < e$.  Put $w = w_{\Bb}(b)$, and
assume that $w \in \s^qW$, (cf. 2.6). Then we have
\begin{equation*}
\Bp_{\Bb, b}(x) = 
  \sum_{(\Ba, \f) \in \wt\CP^q_W}\wt\x^{\Ba,\f}(w)\Bs_{\Ba,\f}(x).
\end{equation*}
\end{prop}
\para{2.12.}
As in [S, 3.5], we consider the ring of symmetric polynomials 
$\Xi_{\Bm} = 
 \bigotimes_{k=0}^{e-1}\ZZ[x_1^{(k)}, \dots,
                     x_{m_k}^{(k)}]^{\FS_{m_k}}$
with respect to $\FS_{\Bm} = \FS_{m_0} \times\cdots\times\FS_{m_{e-1}}$.
$\Xi_{\Bm}$ has a structure of graded ring 
$\Xi_{\Bm} = \bigoplus _{i\ge 0}\Xi_{\Bm}^i$, where $\Xi_{\Bm}^i$
consists of homogeneous symmetric polynomials of degree $i$.  We can 
define a space of symmetric functions $\Xi = \bigoplus_{i\ge 0}\Xi^i$,
where $\Xi^i$ is the inverse limit of $\Xi_{\Bm}^i$ (cf. [loc. cit.]).
For $\Ba \in \CP_{n,e}$, the functions 
$s_{\Ba}(x)$, $p_{\Ba}(x)$ given in 2.2 make sense for 
infinietely many variables 
$x_1^{(k)}, x_2^{(k)}, \dots$,  and give rise to elements in 
$\Xi^n$.  
\par
Schur functions $\Bs_{\Ba,\f}(x)$ and
power sum symmetric functions $\Bp_{\Bb,b}(x)$ defined in 2.10 
are also extended to the functions with infinitely many 
variables, and 
they can be regarded as elements in the space 
$\bigoplus p\Xi_{\CC}^n$, the direct
sum of $p$ copies of $\CC\otimes\Xi^n$.
It is easy to see that Schur functions 
$\{ \Bs_{\Ba,\f} \mid (\Ba,\f) \in \wt\CP^q_W\}$ are linearly
independent in $\bigoplus p\Xi_{\CC}^n$.  We denote by 
$\Xi^n_{\CC}(p,q)$ the
subspace of $\bigoplus p\Xi_{\CC}^n$ generated by those $\Bs_{\Ba,\f}$.
Thus Schur functions $\{\Bs_{\Ba,\f} \mid (\Ba,\f) \in \wt\CP_W^q\}$ 
form a basis of $\Xi^n_{\CC}(p,q)$. 
\par
The space $\Xi^n_{\CC}(p,q)$ is also described as follows.
Let $\Xi^n_{p,q}$ be the subspace of $\CC\otimes\Xi^n$ generated
by $\sum_{i=0}^{c-1}\z^{qid}s_{\th^i(\Ba)}$ for various 
$\Ba \in \CP_{n,e}$ (here $c = c_{\Ba,p}$).  Then
\begin{equation*}
\tag{2.12.1}
\Xi^n_{\CC}(p,q) = \bigoplus_{j=0}^{p-1}\Xi^{n/h_j}_{p,q}(\CX_j), 
\end{equation*}
where $\Xi^{n/h_j}_{p,q}(\CX_j)$ denotes the space $\Xi_{p,q}^{n/h_j}$
with repsect to the variables $\CX_j$. (2.12.1) is immediate from 
(2.10.1). 
\par
It follows from Proposition 2.11 that 
$\Bp_{\Bb,b} \in \Xi^n_{\CC}(p;q)$ 
if $w_{\Bb}(b) \in \s^qW$. By making use of
the Frobenius formula, we can give a complete description of
$W$-orbits of $\s^qW$.
First we prepare some notations.
Take $w = w_{\Bb}(b) \in \s^qW$. For $0 \le j < p$, 
we define $f_j(w) \in \ZZ/p\ZZ$ by 
\begin{equation*}
\tag{2.12.2}
f_j(w) = \begin{cases}
             bj \pmod p  &\quad\text{ if $\Bb$ satisfies (2.6.1)}, \\
             0           &\quad\text{ otherwise}.
         \end{cases}
\end{equation*}  
\par
By Proposition 2.7 and (2.10.2), we can associate to each $W$-orbit
in $\s^qW$ a well-defined function $\Bp_{\Bb,b}(x)$, where 
$w = w_{\Bb}(b)$ is a representative of this orbit. 
We denote $\Bp_{\Bb,b}(x)$ also by $\Bp_w(x)$.
As a corollary to Proposition 2.11, 
we have the following characterization of $W$-orbits in $\s^qW$.
%%%%
\begin{prop}%%%% Prop.2.13.
\begin{enumerate}
\item
Assume that $w, w' \in \s^qW$.
Then $w$ and $w'$ are conjugate under $W$ if and only if 
$\Bp_w(x) = \Bp_{w'}(x)$. 
\item
Let $w = w_{\Bb}(b)$ and $w' = w_{\Bb}(b')$ are elemetns in $\s^qW$, 
(i.e., $\vD(\Bb) \equiv q\pmod p$, cf. 2.3).
Then $w$ and $w'$ are conjugate under 
$W$ if and only if $f_j(w) = f_j(w')$ for $j =0, \dots, p-1$. 
\end{enumerate}
\end{prop}
%%%%
\begin{proof}
It is known that the ``character table'' 
$(\wt\x^{\Ba,\f}(w))$ of $\s^qW$ is a non-singular
matrix. (Here  $(\Ba,\f)$ runs over all the elements
in $\wt\CP_W^q$, and $w$ runs over all the representatives of 
$W$-orbits in  $\s^qW$).  Since Schur functions form a basis of 
$\Xi^n_{\CC}(p,q)$,  we see that $\{\Bp_w \mid w \in \s^qW\ssim\}$ 
is also a basis by Proposition 2.11.  In particular, they are all distinct. 
This proves (i).  The second statement is then immediate from 
(2.10.2).
\end{proof}
%%%%
%%%%
\para{2.14.}
In view of Proposition 2.13, we can determine the parameter set for
the set of $W$-orbits in $\s^qW$.  Put 
\begin{equation*}
\CP_q^W = \{(\Bb,b) \mid \Bb \in \CP_{n,e}, 
              0 \le b <e, \vD(\Bb) \equiv q\pmod p \}.
\end{equation*}
We define an equivalence relation on $\CP_q^W$ by 
$(\Bb,b) \sim (\Bb',b')$ if and only if  
$\Bb = \Bb'$ and $\Bp_{\Bb,b}(x) = \Bp_{\Bb',b'}(x)$.
We denote by 
$\wt\CP_q^W = \wt\CP_q^{n,e,p}$ the
set of equivalnece classes in $\CP_q^W$.  Then $\wt\CP_q^W$
parametrizes the set of $W$-orbits in $\s^qW$, and so the set
$\{ \Bp_{\Bb,b}\mid (\Bb,b) \in \wt\CP_q^W\}$ gives rise to a basis of
$\Xi^n_{\CC}(p,q)$.
\para{2.15.}
In this and next subsection, we give some examples of 
Schur functions and powersum symmetric functions 
in typical cases, i.e., the case where $W = G(2,2,n)$ and the case
where $W = G(e,e,2)$.
\par
First we assume that $W = G(2,2,n)$.  So $W$ (resp. $\wt W$) 
is the Weyl group of type $D_n$ (resp. type $B_n$ ). In this case, 
$\th^2 = 1$, and we have   
$\vG_{\Ba} \simeq \ZZ/2\ZZ$ if $\th(\Ba) = \Ba$, and $\vG_{\Ba} = 1$
otherwise.  In later discussion, the case $\wt\CP^0_W$ corresponds to
the split $D_n$ case, and the case
$\wt\CP_W^1$ corresponds to the non-split $D_n$ case. 
For each $(\Ba,\f) \in \wt\CP^0_W$, Schur function 
$\Bs_{\Ba,\f} = (s^0_{\Ba,\f}, s^1_{\Ba,\f})$ is defined as
follows.
\begin{align*}
\tag{2.15.1}
s^0_{\Ba,\f}(x) &= \sum_{\Bb \in O(\Ba)}s_{\Bb}(x), \\
s^1_{\Ba,\f}(x) &= \begin{cases}
       \pm s_{\Bb}(X) &\quad\text{ if } \th(\Ba) = \Ba, \\
       0              &\quad\text{ otherwise}.
                  \end{cases}   
\end{align*}
Here if $\th(\Ba) = \Ba$, then $\Ba$ is written as $\Ba = (\b;\b)$ with 
a partition $\b$ of $n/2$, and
$s_{\Bb}(X)$ is the Schur function with respect to the variables
$X_j = x_j^{(0)}x_j^{(1)}$. (Note that $\b = \Ba\{1\}$ in the
previous notation).  The sign takes $+$ (resp. $-$) if 
$\f = 1$ (resp. $\f \ne 1$).  
\par
On the other hand, 
$(\Ba,\f) \in \wt\CP^1_W$ if and only if $\th(\Ba) \ne \Ba$, and so
$\vG_{\Ba} = \{1\}$.  In this case, Schur function
$\Bs_{\Ba,\f}(x)$ is given as follows. 
\begin{equation*}
\tag{2.15.2}
s^0_{\Ba,\f}(x) = s_{\Ba}(x) - s_{\th(\Ba)}(x), 
      \qquad s^1_{\Ba,\f}(x) = 0.
\end{equation*} 
\par
Next we consider the power sum symmetric function 
$\Bp_{\Bb,b} = (p^0_{\Bb,b}, p^1_{\Bb,b})$.
A conjugacy class $C_{\Bb}$ of $\wt W$ lying in
$W$ decomposes into two $W$-orbits if and only if 
$\Bb = (\b; -)$ with a partition 
$\b = (\b_1,\b_2,\cdots)$ such that $\b_i$ is even for each
$i$.  Let $\b' = (\b_1/2, \b_2/2, \dots)$ be a partition of $n/2$.
Then we have $\b' = \Bb[1]$ in the previous notation.
Such a class $C_{\Bb}$ is called a degenerate class, and other classes
are said to be non-degenerate.
The function $\Bp_{\Bb,b}(x)$ ($b = 0,1$) associated to 
$(\Bb,b) \in \wt\CP^W_0$ is given as follows.
\begin{align*}
\tag{2.15.3}
p^0_{\Bb,b}(x) &= p_{\Bb}(x), \\
p^1_{\Bb,b}(x) &=  \begin{cases}
         0 &\quad\text{ if $C_{\Bb}$ is non-degenerate,} \\
         \pm 2^{l(\Bb)}p_{\b'}(X)  
            &\quad\text{ if $C_{\Bb}$ is degenerate.} 
                         \end{cases}
\end{align*}  
Here $p_{\Bb}(x)$ is the function for $\wt W$, and 
$p_{\b'}(X)$ is the powersum symmetric function for 
$\FS_{n/2}$ with variables 
$X_j = x_j^{(0)}x_j^{(1)}$ associated to the partition $\b'$.
The sign takes $+$ (resp. $-$) if $b = 0$ (resp. $b = 1$).
\par
On the other hand, $(\Bb, b) \in \wt\CP^W_1$ if and only if 
$\vD(\Bb) \equiv 1 \pmod 2$ and $b = 0$.  In this case, 
$\Bp_{\Bb,b}(x)$ is given as
\begin{equation*}
p_{\Bb,b}^0(x) = p_{\Bb}(x), \qquad p_{\Bb,b}^1(x) = 0.
\end{equation*}
%%%%
\para{2.16.}
We now consider the case where $W = G(e,e,2)$.  So $W$ is 
the dihedral group of order $2e$.  
Put 
\begin{align*}
\Ba_0 &= (2;-;\cdots;-), \\ 
\Ba_{ij} &= (-;\cdots;1;\cdots;1;\cdots;-) \quad (0\le i\le j <e),
\end{align*}
where in the second case, 1 appears only in the $i$-th and $j$-th
entries. If $i = j$, we understand that 
$\Ba_{ii} = (-;\cdots;11;\cdots;-)$.
Then  the set of 
representatives of $\th$-orbits of $\CP_{2,e}$ is given as 
$\{ \Ba_0, \Ba_{0j} \mid 0 \le j \le e/2 \}$.  
Here $\vG_{\Ba} = \{1\}$ unless $e$ is even and
$j = e/2$, in which case $\vG_{\Ba_{0j}} = \ZZ/2\ZZ$.
Consequently, Schur functions 
$\Bs_{\Ba,\f} = (s^k_{\Ba,\f})_{0 \le k <e}$ associated to
$(\Ba,\f) \in \wt\CP^0_W$ are given as follows.
\begin{equation*}
\tag{2.16.1}
s^k_{\Ba,\f}(x) = \begin{cases}
      \sum_{\Bb \in O(\Ba)}s_{\Bb}(x) &\quad\text{ if } k = 0, \\
      \pm\sum_{\Bb' \in O(\Ba')}s_{\Bb'}(X)
           &\quad\text{ if $\Ba = \Ba_{0,e/2}$ and $k = e/2$}, \\
       0   &\quad\text{ otherwise.}
                   \end{cases}
\end{equation*} 
In the second case, $\Ba' = (1;-;\cdots;-)$ is an $e/2$-partition of $1$ and
coincides with $\Ba\{e/2\}$. $s_{\Ba'}(X)$ is the Schur function
associated to $\Ba'$ with respect to the variables 
$X^{(0)},\dots, X^{(e/2-1)}$ for $X_j^{(k)} = x_j^{(k)}x_j^{(e/2+k)}$.
Hence $s_{\Ba'}(X) = \sum_j x_j^{(0)}x_j^{(e/2)}$ and we have
\begin{equation*}
s^{e/2}_{\Ba,\f}(x) = \pm\sum_{i=0}^{e/2-1}\sum_jx_j^{(i)}x_j^{(i+e/2)}.
\end{equation*}  
The sign takes $+$ (resp. $-$) if $\f = 1$ (resp. $\f \ne 1$).
\par
Next we consider power sum symmetric functions. 
A conjugacy class $C_{\Bb} \subset \wt W$ belongs to $W$ in the
following case; 
\begin{equation*}
\{ \Ba_0, \quad \Ba_{i,e-i} \ (0 \le i \le e/2) \}.
\end{equation*}
\par
Then $\Bp_{\Bb,b} = (p^k_{\Bb,b}(x))_{0 \le k < e}$ 
associated to $(\Bb,b) \in \wt\CP^W_0$ is given as
\begin{equation*}
\tag{2.16.2}
p^k_{\Bb,b}(x) = \begin{cases}
                  p_{\Bb}(x) &\quad\text{ if } k = 0, \\
                  (-1)^b2p_{\Ba'}(X)
          &\quad\text{ if  $k = e/2$ and  $\Bb = \Ba_0$},  \\
                   0  &\quad\text{ otherwise},
                  \end{cases}   
\end{equation*} 
where $\Ba' = (1;-;\cdots;-)$  is 
the $e/2$-partition as given before and also coincides with 
$\Ba_0[e/2]$, 
and
$p_{\Ba'}(X)$ is the power sum symmetric function associated to $\Ba'$
with respect to variables $\{ X^{(0)}, \dots, X^{(e/2-1)}\}$.  
Hence we have 
\begin{equation*}
p_{\Ba'}(X) = \sum_{i=0}^{e/2-1}\sum_jx_j^{(i)}x_j^{(i+e/2)}.
\end{equation*}
\par\medskip
\section{Hall-Littlewood functions associated to $G(e,p,n)$}
\para{3.1.}
Here we review some results from [S]  
concerning symmetric functions associated to $\wt W$.
In what follows we regard
the variables $x_i^{(k)}$ defined for 
$k \in \ZZ/e\ZZ \simeq \{ 0,1, \dots, e-1\}$.  
For each $0 \le k < e$ 
and an integer $r \ge 0$, 
we define a function $q_{r,\pm}^{(k)}(x;t)$ (according to the sign $+$
or $-$) by 
\begin{equation*}
\tag{3.1.1}
q_{r,\pm}^{(k)}(x;t) = \sum_{i\ge 1} (x_i^{(k)})^{r+\d}
       \frac{\prod_{j}x_i^{(k)} -tx_j^{(k\pm 1)}}
     {\prod_{j\ne i}x_i^{(k)} - x_j^{(k)}}  \qquad (r \ge 1),
\end{equation*} 
where $\d = m_{k} - 1 - m_{k\pm 1}$,
and by $q_{r,\pm}^{(k)}(x;t) = 1$ for $r = 0$. In the product of 
the denominator,
$x_j^{(k)}$ runs over all the variables in $x^{(k)}$ except 
$x_i^{(k)}$, while in the numerator,
$x_j^{(k\pm 1)}$ runs over all the variables in $x^{(k\pm 1)}$.
Then $q_{r,\pm}^{(k)}(x;t)$ is a polynomial in $\ZZ[x;t]$,  
homogenoeous of degree $r$ with respect to the variables 
$x^{(k)}, x^{(k\pm 1)}$ ([S, Lemma 2.3]).
\par
For an $e$-partition 
$\Ba = (\a^{(0)}, \dots, \a^{(e-1)}) \in \CP_{n,e}$, 
we define a function $q_{\Ba,\pm}(x)$ by
\begin{equation*}
\tag{3.1.2}
q_{\Ba,\pm}(x;t) = \prod_{k=0}^{e-1}
        \prod_{j=1}^{m_k}q_{\a_j^{(k)},\pm}^{(k)}(x;t).
\end{equation*}
\par
For $\Ba = (\a_j^{(k)}) \in \CP_{n.e}$, we define a function 
$z_{\Ba}(t)$ by 
\begin{equation*}
\tag{3.1.3}
z_{\Ba}(t) = 
z_{\Ba}\prod_{k=0}^{e-1}\prod_{j=1}^{m_k}(1 - \z^kt^{\a^{(k)}_j})\iv,
\end{equation*}
where in the product, we neglect the factors such that 
$\a_j^{(k)} = 0$. $z_{\Ba}$ is the order of the centralizer 
of $w_{\Ba}$ in $\wt W$.
Explicitly, $z_{\Ba}$ is given as follows.  For a partition 
$\b = (1^{n_1}, 2^{n_2}, \dots)$, put  
$z_{\b} = \prod_{i\ge 1}i^{n_i}n_i!$.  Then 
$z_{\Ba} = e^{l(\Ba)}\prod_{k=0}^{e-1}z_{\a^{(k)}}$.
\par
For later use, we also give the
order of the stabilzer $Z_W(w)$ of $w \in \s^qW$ in $W$.  
Assume that $C_{\Ba} \subset \s^qW$ and that $C_{\Ba}$ decomposes 
into $r$ distinct $W$-orbits.  Since these $W$-orbits have the same
cardinality, we have for $w \in C_{\Ba}$,
\begin{equation*}
\tag{3.1.4}
|Z_W(w)| = \frac {r}{p}|Z_{\wt W}(w)|.
\end{equation*}
\par
We now introduce infinitely many variables $x_i^{(k)}, y_i^{(k)}$ for
$i = 1,2, \dots$ and for $0 \le k \le e-1$.  As discussed in [S], 
$p_{\Ba}(x), q_{\Ba,\pm}(x;t), m_{\Ba}(x)$ can be regarded as
functions  with infinitely
many variables $x_1^{(k)}, x_2^{(k)}, \dots$.  
Following [S], we introduce  Cauchy's reproducing kernel 
associated to $\wt W$ by
\begin{equation*}
\tag{3.1.5}
\Om(x,y;t) = \prod_{k=0}^{e-1}\prod_{i,j}
    \frac{1 - tx_i^{(k+1)}y_j^{(k)}}{1 - x_i^{(k)}y_j^{(k)}}.
\end{equation*}
The following formula was proved in [S, Proposition 2.5].
%%%
\begin{prop}%%%Prop. 3.2
We have
\begin{align*}
\tag{3.2.1}
\Om(x,y;t) &= \sum_{\Ba}q_{\Ba,+}(x;t)m_{\Ba}(y) 
                = \sum_{\Ba}m_{\Ba}(x)q_{\Ba,-}(y;t), \\
\tag{3.2.2}
\Om(x,y;t) &= \sum_{\Ba}z_{\Ba}(t)\iv p_{\Ba}(x)\bar p_{\Ba}(y),
\end{align*}
where $\Ba$ runs over all the $e$-partitions 
of any size. In (3.2.2), $\bar p_{\Ba}(y)$ denotes the complex 
conjugate of $p_{\Ba}(y)$.
\end{prop}   
\para{3.3.}
We define an isomorphism $\th$ of $\CC[x]$ (the polynomial ring with
infinitely many variables) by $\th(x_j^{(k)}) = x_j^{(k+d)}$.
Then it follows from (2.2.2) that 
$\th(p_r^{(k)}) = \z^{-kd}p_r^{(k)}$, and so we have
\begin{equation*}
\tag{3.3.1}
\th(p_{\Ba}) = \z^{-\vD(\Ba)d}p_{\Ba}.
\end{equation*}
\par 
On the other hand, it can be checked easily from the definition that
\begin{equation*}
\tag{3.3.2}
\th(s_{\Ba}) = s_{\th(\Ba)}, \qquad \th(m_{\Ba}) = m_{\th(\Ba)}, 
\qquad \th(q_{\Ba,\pm}) = q_{\th(\Ba),\pm}.
\end{equation*}
\par
We now introduce an operator $\vT_q: \CC[x] \to \CC[x]$ for 
a factor $q$ of $e$ by
\begin{equation*}
\tag{3.3.3}
\vT_q = \sum_{i = 0}^{p-1}\z^{qid}\th^i.
\end{equation*}
Then it is easy to see that 
\begin{equation*}
\tag{3.3.4}
\vT_q(p_{\Bb}) = \begin{cases}  
                      p\cdot p_{\Bb}  &\quad\text{ if } 
                                 \vD(\Bb) \equiv q \pmod p, \\
                      0        &\quad\text{ otherwise}.
                  \end{cases}   
\end{equation*}
Take $\Ba \in \CP_{n,e}$ and put $c = c_{\Ba,p}$ as before. 
Then it follows from (3.3.2) that we have
\begin{equation*}
\tag{3.3.5}
\vT_q(s_{\Ba}) = \begin{cases} 
       \displaystyle{\frac p c\sum_{i=0}^{c-1}\z^{qid}s_{\th^i(\Ba)}}
            &\quad\text{ if } qc \equiv 0\pmod p, \\
       0    &\quad\text{ otherwise}.
                  \end{cases} 
\end{equation*}
Similar formulas as (3.3.5) hold also for $\vT_q(m_{\Ba})$ and 
$\vT_q(q_{\Ba,\pm})$ by replacing $s_{\th^i(\Ba)}$ in (3.3.5) by
$m_{\th^i(\Ba)}$ and $q_{\th^i(\Ba),\pm}$.
\par
Let us define $m^j_{\Ba,\f}(x)$ and
$q^j_{\Ba, \f,\pm}(x;t)$, for each $0 \le j < p$,  as in the case
of $\Bs_{\Ba,\f}(x)$ by replacing $s_{\th^i(\Ba)\{j\}}$ 
in (2.10.1) by $m_{\th^i(\Ba)\{j\}}$ and  $q_{\th^i(\Ba)\{j\},\pm}$. 
We put 
\begin{equation*}
\tag{3.3.6}
\Bq_{\Ba,\f,\pm} = (q^j_{\Ba,\f, \pm})_{0 \le j < p}, \qquad
\Bm_{\Ba,\f} = (m^j_{\Ba,\f})_{0 \le j < p}.
\end{equation*}
Note that $\th$ induces an action  
$\th(X_i^{(k)}) = X_i^{(k+d)}$ on the variables $\CX_j$.
Thus similar formulas as (3.3.4) and (3.3.5) hold also for funlctions
$p_{\Ba[j]}(\CX_j)$ and $s_{\Ba\{j\}}(\CX_j)$, etc. 
(More precisely, the formula for $s_{\Ba\{ j\}}$ is completely similar
to (3.3.5).  But for $p_{\Bb}(\CX_j)$, $\z$ must be replaced by
$\z^h$.  In particular, we have 
$\th(p_r^{(k)}(\CX_j)) = \z^{-khd}p_r^{(k)}(\CX_j)$.
Hence the condition in the first formula in (3.3.4) should be
replaced by ``if $h\vD(\Bb) \equiv q \pmod p$'' for a 
$j_1d$-partition $\Bb$).
Now $q_{\Ba,\pm}$ can be written as a linear combination of
$s_{\Bb}$.  Then by applying $\vT_q$, 
$\sum_i\z^{qid}q_{\th^i(\Ba),\pm}$ can be written as a linear 
combination  of various $\sum_i\z^{qid}s_{\th^i(\Bb)}$. 
It follows that 
$q_{\Ba,\f,\pm}^j \in \Xi_{p,q}^{n/h_j}(\CX_j)$.    
Hence by (2.12.1), $\Bq_{\Ba,\f,\pm}$ is 
contained in $\Xi^n_{\CC(t)}(p,q) = \CC(t)\otimes \Xi^n_{\CC}(p,q)$. 
 Similarly, $\Bm_{\Ba,\f}$ 
is contained in $\Xi^n_{\CC(t)}(p,q)$. Now it is easy to see that 
the sets
$\{ \Bq_{\Ba,\f,\pm} \mid (\Ba,\f) \in \wt\CP^q_W\}$ and
$\{ \Bm_{\Ba,\f} \mid (\Ba,\f) \in \wt\CP^q_W\}$ form
bases of $\Xi^n_{\CC(t)}(p,q)$.
\par
For each $j$ ($0 \le j < p)$, we also consider the variables 
$\CY_j = (Y_i^{(k)})$ defined in a similar way as $\CX_j$ but 
obtained from the variables $y = (y_i^{(k)})$ instead of 
$x = (x_i^{(k)})$.  
We now define Cauchy's reproducing kernel associated to $\s^qW$ by
$\BOm_q(x,y;t) = (\Om_q^j(x,y;t))_{0 \le j <p}$, with
\begin{equation*}
\tag{3.3.7}
\Om_q^j(x,y;t) = \vT_{q,\CX_j}(\Om(\CX_j,\CY_j;t^{h_j})),  
\end{equation*} 
where $\Om(\CX_j,\CY_j;t^{h_j})$ is the function defined in a similar way 
as (3.1.4), with respect to the variables $\CX_j, \CY_j$ and $t^{h_j}$. 
$\vT_{q,\CX_j}$ stands for the action of $\vT_q$ on the variables $\CX_j$.
For $w = w_{\Ba}(b) \in \s^qW$, we define $z_{\Ba,b}(t)$ by 
\begin{equation*}
\tag{3.3.8}
z_{\Ba,b}(t) = 
  z_{\Ba,b}\prod_{k=0}^{e-1}\prod_{j=1}^{m_k}(1-\z^kt^{\a_j^{(k)}})\iv.
\end{equation*}
Hence, $z_{\Ba,b}(t)$ is the one obtained from $z_{\Ba}(t)$ by
replacing  $z_{\Ba} = |Z_{\wt W}(w)|$ by $z_{\Ba,b} = |Z_W(w)|$, and 
so it coincides with $rp\iv Z_{\Ba}(t)$ by (3.1.4).
\par
The following proposition gives a counter part of 
Propostion 3.2 to the case of $\s^qW$.
%%%%
\begin{prop} %%%% Prop. 3.4
$\BOm_q(x,y;t)$ has the following expansions.
\begin{align*}
\tag{3.4.1}
\BOm_q(x,y;t) &= \sum_{(\Ba,\f)}\Bq_{\Ba,\f,+}(x;t)\ol{\Bm}_{\Ba,\f}(y) 
                = \sum_{(\Ba,\f)}\Bm_{\Ba,\f}(x)\ol{\Bq}_{\Ba,\f,-}(y;t), \\
\tag{3.4.2}
\BOm_q(x,y;t) &= \sum_{(\Ba,b)}z_{\Ba,b}(t)\iv \Bp_{\Ba,b}(x)\ol{\Bp}_{\Ba,b}(y),
\end{align*}
where $(\Ba,\f)$ runs over all the elements in 
$\bigcup_{n=1}^{\infty} \wt\CP^q_{n,e,p}$ in (3.4.1). 
In (3.4.2), $(\Ba,b)$ runs over all the elements in
$\bigcup_{n=1}^{\infty}\wt\CP^{n,e,p}_q$.
$\ol{\Bm}_{\Ba,\f}$ (resp. $\ol{\Bq}_{\Ba,\f,-}$,
$\ol{\Bp}_{\Ba,b}$)
 denotes the complex conjugate of $\Bm_{\Ba,\f}$ 
(resp. $\Bq_{\Ba,\f,-}$, $\Bp_{\Ba,b}$), respectively.
\end{prop}
%%%%  
\begin{proof}
First we show (3.4.1).
We fix $j$ and put $h = h_j$ as before.  Let us consider the expansion 
of $\Om(\CX_j,\CY_j;t^h)$ by 
making use of
the first equality of (3.2.1). By applying the operator $\vT_{q,\CX_j}$ 
on both sides of this expansion, we have
\begin{equation*}
\tag{3.4.3}
\Om_q^j(x,y;t) = \sum_{\Bb}\vT_q(q_{\Bb,+}(\CX_j;t^h))m_{\Bb}(\CY_j),
\end{equation*}
where $\Bb$ runs over $j_1d$-partitions of any size.
Take $\Bb \in \CP_{n/h, j_1d}$.
In view of (3.3.5), we may only consider
$\Bb$ such that $qc \equiv 0 \pmod {p}$ in the sum (3.4.3), where
$c = c_{\Bb,p}$.  
By making use of the formula for $\vT_q(q_{\Bb,+})$ which is similar
to (3.3.5), we have
\begin{align*}
\tag{3.4.4}
\sum_{\Bb' \in O(\Bb)}\vT_q(q_{\Bb',+})m_{\Bb'}   
    &= \sum_{\Bb' \in O(\Bb)}\biggl\{\frac p {c}\sum_{i=0}^{c-1}\z^{qid}
             q_{\th^i(\Bb'),+}\biggr\}m_{\Bb'} \\
     &= \frac p {c}\sum_{i=0}^{c-1}\z^{qid}q_{\th^i(\Bb),+}
                 \sum_{k=0}^{c-1}\z^{-qkd}m_{\th^k(\Bb)}.
\end{align*}  
Now there exists $\Ba \in \CP_{n,e}$ such that $\Bb = \Ba\{j\}$.  Then
$c$ coinicdes with $c_{\Ba,p}$. Moreover the action of $\th$ on $\Bb$
is compatible with that on $\Ba$.
Then the last formula of (3.4.4) can be
written as 
\begin{align*} 
&\sum_{\f \in \vG_{\Ba}\wg}
   \biggl\{\f(\t^j)\sum_{i=0}^{c-1}
     \z^{qid}q_{\th^i(\Ba)\{j\},+}(\CX_j;t^h)\biggr\}
 \biggl\{\f(\t^{-j})\sum_{k=0}^{c-1}
        \z^{-kid}m_{\th^i(\Ba)\{j\}}(\CY_j)\biggr\}   \\
= &\sum_{\f \in \vG_{\Ba}\wg}q^j_{\Ba,\f,+}(x;t))\ol {m^j_{\Ba,\f}(y)}.   
\end{align*}
It follows that
\begin{equation*}
\Om_q^j(x,y;t) = \sum_{(\Ba,\f)}
        q^j_{\Ba,\f,+}(x;t)\ol{m^j_{\Ba,\f}(y)},
\end{equation*}
where $(\Ba,\f)$ runs over all the elements in 
$\bigcup_{n=1}^{\infty}\wt\CP^q_{n,e,p}$.
This shows the first equality of (3.4.1).  The second equality 
is shown similarly. 
\par
Next we show (3.4.2).  We consider the expansion of 
$\Om(\CX_j,\CY_j;t^h)$ by (3.2.2).  By applying $\vT_{q,\CX_j}$ on this
equality, together with (3.3.4), we have
\begin{equation*}
\tag{3.4.5}
\Om_q^j(x,y;t) = p\sum_{\Bb}
    z_{\Bb}(t^h)\iv p_{\Bb}(\CX_j)\ol{p_{\Bb}(\CY_j)}, 
\end{equation*}
where $\Bb$ runs over $j_1d$-partitions of any size such that
$h\vD(\Bb) \equiv q \pmod p$ (see the discussion in 3.3).  
$z_{\Bb}(t^h)$ is the function defined similar to (2.1.3) 
for a $j_1d$-partition $\Bb$, by replacing $\z$ by $\z^h$. 
We fix $\Bb \in \CP_{n/h, j_1d}$ as above.  Then there 
exists $\Ba \in \CP_{n,e}$ such
that $\Bb = \Ba[j]$.  Hence $\Ba$ satsifies (2.6.1), and 
$\vD(\Ba) = h\vD(\Bb)$.  Let $r$  be the number of all the pairs 
$(\Ba, b) \in \wt\CP^{n,e,p}_q$ for a fixed $\Ba$.   
Then by (2.10.2), we have
\begin{equation*}
\tag{3.4.6}
p_{\Bb}(\CX_j)\ol{p_{\Bb}(\CY_j)} = 
h^{-2l(\Bb)}r\iv\sum_b p^j_{\Ba,b}(x)\ol{p^j_{\Ba,b}(y)},
\end{equation*}
where the sum is taken over all the pairs $(\Ba,b)$ for a fixed $\Ba$.
Now using the explicit description of $z_{\Ba}(t)$ in (3.1.3) and
subsequent parts, one can check that
\begin{equation*}
\tag{3.4.7}
z_{\Bb}(t^h) = h^{-2l(\Bb)}z_{\Ba}(t).
\end{equation*}  
Substituting (3.4.6) and (3.4.7) into (3.4.5), together with (3.3.8),
we have
\begin{equation*}
\Om_q^j(x,y;t) = \sum_{(\Ba,b)}
            z_{\Ba,b}(t)\iv p^j_{\Ba,b}(x)\ol{p^j_{\Ba,b}(y)},
\end{equation*}
where $(\Ba,b)$ runs over all the elements in 
$\bigcup_{n=1}^{\infty}\wt\CP_q^{n,e,p}$.  
This implies (3.4.2), and the propostion is proved.
\end{proof}
%%%%
\para{3.5.}
As in [S], we denote by 
$Z_n^{0,0} = Z_n^{0,0}(\Bm)$   
the set of e-partitions $\Ba$ such that $|\Ba| = n$ and that
each $\a^{(k)}$ is regarded as an element in $\ZZ^{m_k}$, 
written in the form 
$\a^{(k)}: \a^{(k)}_1 \ge \cdots \ge \a^{(k)}_{m_k} \ge 0$.       
We fix an integer $r > 0$, and consider an 
$e$-partition $\BLa^0 = \BLa^0(\Bm) = (\vL_0, \dots, \vL_{e-1})$
defined by 
\begin{equation*}
\tag{3.5.1}
\vL_i: \ (m_i-1)r > \cdots >  2r > r > 0
                    \qquad \text{ for } 0 \le i \le e-1. 
\end{equation*}
In [S], the set $Z_n^{r,s}$ of symbols are introduced.  
In this paper, we are only concerned with the 
symbol of the form $Z_n^{r,0} = Z_n^{r,0}(\Bm)$, i.e., 
 the set of $e$-partitions of the form 
$\BLa = \Ba + \BLa^0$, where $\Ba \in Z_n^{0,0}$ and the 
sum is taken entry-wise.
We denote by $\BLa = \BLa(\Ba)$ if $\BLa = \Ba + \BLa^0$, and 
call it the $e$-symbol of type $(r,0)$ corresponding 
to $\Ba$.  We write $|\BLa| = n$ if $\BLa \in Z_n^{r,0}$.
\par
As in [S], we consider the shift operation
$Z_n^{r,0}(\Bm) \to Z_n^{r,0}(\Bm')$ for
$\Bm' = (m_0+1, \dots, m_{e-1}+1)$ by 
associating $\BLa' = \Ba + \BLa^0(\Bm')$
to $\BLa = \Ba + \BLa^0(\Bm)$, where $\Ba$ is 
regarded as an element of $Z_n^{0,0}(\Bm')$ by adding 0 in the entries
of $\Ba$. We often consider the symbols as elements in the equivalence
class in the set $\coprod_{\Bm'}Z_n^{r,0}(\Bm')$ under the shift 
operation. We denote by $\bar Z_n^{r,0}$ the set of equivalence
classes.   
Note that $\BLa^0$ is regarded as a symbol in $Z_n^{r,0}$ with 
$n = 0$.
\par
Two elements $\BLa$ and $\BLa'$ in $\bar Z_n^{r,0}$ are said to be similar
(and denoted as $\BLa \sim \BLa'$) if there exist representatives 
in $Z_n^{r,0}$ such that 
all the entries of them coincide each other with multiplicities. 
The equivalence class with respect to this relation is called a
similarity class in $\bar Z_n^{r,0}$.
\par
We shall define a function $a : \bar Z_n^{r,0} \to \NN$.  For 
$\BLa \in  Z_n^{r,0}$, we put 
\begin{equation*}
\tag{3.5.2}
a(\BLa) = \sum_{\la,\la' \in \BLa}\min(\la,\la')
            - \sum_{\m,\m' \in \BLa^0}\min(\m,\m'),
\end{equation*}
which induces a well-defined function $a$ on  
$\bar Z_n^{r,0}$. 
The $a$-function takes a constant value on each similarity
class in $Z_n^{r,0}$.
By using the bijection 
$Z_n^{0,0} \simeq Z_n^{r,0}$, 
$a$-functions and similarity classes are defined for $Z_n^{0,0}$, for
which we use the same notation as for $Z_n^{0,0}$.
\par
Under the natural bijection 
$\CP_{n,e} \simeq Z_n^{0,0} \simeq Z_n^{r,0}$, the operation $\th$ is
transfered to the action on $Z_n^{r,0}$, which we denote also by
$\th$.  It is clear that $\th$ preserves each similarity class 
in $Z_n^{r,0}$.
\para{3.6.}
As discussed in [S], $\{ s_{\Ba}(x)\}$ and $\{ m_{\Ba}(x)\}$ form  
bases of the $\ZZ[t]$-module $\ZZ[t] \otimes \Xi$, and 
$\{ q_{\Ba,\pm}(x;t)\}$
form a basis of $\QQ(t) \otimes _{\ZZ}\Xi = \Xi_{\QQ}[t]$.  
Moreover, $\{ p_{\Ba}(x)\}$ gives a
basis of $\CC(t)$-space $\Xi_{\CC}[t] = \CC(t)\otimes_{\ZZ}\Xi$.
\par
Following [S], we define a scalar product on 
$\Xi_{\QQ}[t]$ by the condition that 
\begin{equation*}
\tag{3.6.1}
\lp q_{\Ba, +}(x;t), m_{\Bb}(x)\rp = \d_{\Ba,\Bb},
\end{equation*}
and extend it to a sesquilinear form on $\Xi_{\CC}[t]$.
Then by Proposition 3.2, we also have 
\begin{align*}
\tag{3.6.2}
\lp m_{\Ba}(x), q_{\Bb,-}(x;t)\rp &= \d_{\Ba,\Bb}, \\
\lp p_{\Ba}(x), p_{\Bb}(x) \rp &= z_{\Ba}(t)\d_{\Ba, \Bb}.
\end{align*}
\par
The Hall-Littlewood functions $P_{\BLa}^{\pm}(x;t)$ and 
$Q_{\BLa}^{\pm}(x;t)$ were introduced
in [S].  
By Corollary 4.6 in [S], they satisfy the following formulas.
\begin{align*}
\tag{3.6.3}
\Om(x,y;t) &= \sum_{\BLa, \BLa'}b_{\BLa,\BLa'}(t)
                        P_{\BLa}^+(x;t)P_{\BLa'}^-(y;t). \\
\tag{3.6.4}
\Om(x,y;t) &= \sum_{\BLa}Q^+_{\BLa}(x;t)P^-_{\BLa}(y;t)
           = \sum_{\BLa}P_{\BLa}^+(x;t)Q_{\BLa}^-(y;t),
\end{align*}
where in $(3.6.3)$, $\BLa, \BLa'$ run over all the elements in 
$\bigcup_{n = 1}^{\infty}Z_n^{r,0}$, and $b_{\BLa,\BLa'}(t) = 0$
unless $|\BLa| = |\BLa'|$ and $\BLa \sim \BLa'$.  In $(3.6.4)$, 
$\BLa$ runs over all the elements in 
$\bigcup_{n = 1}^{\infty}Z_n^{r,0}$.
\par
We now define a total order $\prec$ on the set $Z_n^{0,0}$ satisfying
the condition that each similarity class forms an interval and that 
$a(\Ba) > a(\Bb)$ implies that $\Ba \prec \Bb$.
The following characterization of Hall-Littlewood functions 
was given in Propostion 4.8 in [S].
\begin{prop}%%% Prop. 3.7.
\begin{enumerate}
\item
$P_{\BLa}^{\pm}$ are characterized by the following two properties.
\begin{enumerate}
\item
$P^{\pm}_{\BLa}(x;t)$ can be expressed 
in terms of $s_{\Bb}(x)$ as
\begin{equation*}
P^{\pm}_{\BLa} = s_{\Ba} + \sum_{\Bb}u^{\pm}_{\Ba, \Bb}s_{\Bb} 
\end{equation*}
with $u^{\pm}_{\Ba,\Bb} \in \QQ(t)$, where 
$\BLa = \BLa(\Ba)$, and $u^{\pm}_{\Ba,\Bb} = 0$ unless 
$\Bb \prec \Ba$ and $\Bb \not\sim \Ba$.
\item
$\lp P^{+}_{\BLa}, P^-_{\BLa'}\rp = 0$ unless $\BLa \sim \BLa'$.
\end{enumerate}
\item
The functions $\{Q_{\BLa}^{\pm}\}$ are characterized as the dual
basis of $\{ P_{\BLa}^{\pm}\}$, i.e., 
\begin{equation*}
\lp P_{\BLa}^+, Q_{\BLa'}^-\rp = \lp Q_{\BLa}^+, P_{\BLa'}^-\rp
                               = \d_{\BLa,\BLa'}.
\end{equation*}
\end{enumerate}
\end{prop}
\para{3.8.}
We now define a sesquilinear form on 
$\Xi^n_{\CC(t)}(p,q)$ by the condition that 
\begin{equation*}
\tag{3.8.1}
\lp \Bq_{z, +}(x;t), \Bm_{z'}(x)\rp = \d_{z,z'}
\end{equation*}
for $z, z' \in \wt\CP^q_{n,e,p}$.
Then by Proposition 3.4, we also have 
\begin{align*}
\tag{3.8.2}
\lp \Bm_z(x), \Bq_{z'}(x;t)\rp &= \d_{z, z'} 
     \qquad (z, z' \in \wt\CP^q_{n,e,p}),  \\
\lp \Bp_{\xi}(x), \Bp_{\xi'}(x) \rp &= z_{\xi}(t)\d_{\xi, \xi'} 
     \qquad (\xi, \xi' \in \wt\CP_q^{n,e,p}).   
\end{align*}
Let us denote by $\wt Z^{0,0}_{n,q}$ (resp. $\wt Z_{n,q}^{r,0}$ )
the set of pairs $(\Ba, \f)$ (resp. $(\BLa(\Ba),\f)$) 
such that $(\Ba,\f) \in \wt\CP_{n,e,p}^q$ 
and that $\Ba \in Z_n^{0,0}$.  We often identify 
$\wt Z^{0,0}_{n,q}$ and $\wt Z_{n,q}^{r,0}$ with $\wt\CP_{n,e,p}^q$.
A similarity class and an $a$-function on $\wt Z_{n,q}^{0,0}$
are defiend as those inherited from $Z_n^{0,0}$.  
We define a total order $\prec$
on $\wt Z_{n,q}^{0,0}$ compatible with the order on $Z_n^{0,0}$,
i.e., for $z = (\Ba,\f), z' = (\Ba', \f') \in \wt Z_{n,q}^{0,0}$, we have
$z \prec z'$ if $a(\Ba) > a(\Ba')$, and the set
$\{ (\Ba,\f) \}$, where $\Ba$ falls in a fixed similarlity class in 
$Z_n^{0,0}$, form an interval in this order.  
\par
The scalar product on $\Xi_{\CC}[t]$ given in 3.6 induces a 
sesquilinear form $\lp\ ,\ \rp_j$ on the space 
$\CC(t) \otimes\Xi_{p,q}^{n/h_j}(\CX_j)$, which 
satisfies the formula
\begin{equation*}
\tag{3.8.3}
\lp q^j_{\Ba,\f,+}(x), m^j_{\Ba',\f'}(x)\rp_j =
     \begin{cases}
   \f(\t^j)\ol{\f'(\t^{j})}c_{\Ba,p} &\text{ if $\Ba = \Ba'$ 
                  and $c_{\Ba,p} \mid j$}, \\
     0                   &\text{ otherwise}.
    \end{cases} 
\end{equation*}
This implies that for 
$\Bu = (u^j), \Bv = (v^j) \in 
\Xi^n_{\CC(t)}(p,q) = 
  \bigoplus_{j=0}^{p-1}\CC(t)\otimes\Xi^{n/h_j}_{p,q}(\CX_j)$, 
the following formula holds, since
it certainly holds for $\Bu = \Bq_{\Ba,\f,+}$ and
$\Bv = \Bm_{\Ba',\f'}$ by (3.8.3).
\begin{equation*}
\tag{3.8.4}
\lp \Bu, \Bv\rp = 
     \frac 1 p\sum_{j=0}^{p-1}\lp u^j, v^j \rp_j. 
\end{equation*}
\para{3.9.}
In view of (3.3.2) and (3.6.1), the automorphism 
$\th$ leaves the scalar product on $\Xi_{\QQ}[t]$ invariant.  
It follows by Proposition 3.7 that we have
\begin{equation*}
\th(P_{\BLa}^{\pm}) = P_{\th(\BLa)}^{\pm}, \qquad
\th(Q_{\BLa}^{\pm}) = Q_{\th(\BLa)}^{\pm}.
\end{equation*} 
Hence a similar formula as (3.3.5) holds also for
$P_{\BLa}^{\pm}$ (resp. for $Q_{\BLa}^{\pm}$), by replacing 
$s_{\th^i(\Ba)}$ by $P^{\pm}_{\th^i(\BLa)}$ (resp. by 
$Q^{\pm}_{\th^i(\BLa)}$).
Now for each $z = (\Ba,\f)\in \wt Z_{n,q}^{0,0}$, we define 
\begin{equation*}
\tag{3.9.1}
\BP_z^{\pm} = (P_z^{\pm,j})_{0 \le j <p}, \qquad
\BQ_z^{\pm} = (Q_z^{\pm.j})_{0 \le j < p},
\end{equation*}
where $P_z^{\pm,j}$ and
$Q_z^{\pm,j}$ are defined as in the case
of $\Bs_{\Ba,\f}(x)$ by replacing $s_{\th^i(\Ba)\{j\}}$ 
in (2.10.1) by $P^{\pm}_{\th^i(\BLa)\{j\}}$ and  
$Q^{\pm}_{\th^i(\BLa)\{j\}}$. 
By using a s similar argument as in 3.3, we see that
the sets $\{ \BP_z^{\pm}\}$ and $\{ \BQ_z^{\pm}\}$
give rise to bases of 
$\Xi_{\CC(t)}^n(p,q)$. 
\par
It is known that $P^{\pm}_{\BLa}(x;0) = s_{\Ba}(x)$ for 
$\BLa(\Ba) = \BLa$ by [S].  Hence by (3.9.1), we have
\begin{equation*}
\tag{3.9.2}
\BP^{\pm}_{\BLa(\Ba),\f}(x;0) = \Bs_{\Ba,\f}(x). 
\end{equation*}
\par
The functions $\BP^{\pm}_z(x)$ and $\BQ^{\pm}_z(x)$
are charactersized by the following properties, in analogy to 
Proposition 3.7.
\begin{prop}%%% Prop. 3.10.
\begin{enumerate}
\item
$\BP_{z}^{\pm} \ (z \in \wt Z^{r,0}_{n,q})$ 
are characterized by the following two properties.
\begin{enumerate}
\item
$\BP^{\pm}_{z}(x;t)$ can be expressed 
in terms of $\Bs_{z'}(x)$ as
\begin{equation*}
\BP^{\pm}_{z} = \Bs_{z} + \sum_{z'}u^{\pm}_{z, z'}\Bs_{z'} 
\end{equation*}
with $u^{\pm}_{z, z'} \in \CC(t)$, 
and $u^{\pm}_{z,z'} = 0$ unless $z' \prec z$ and $z' \not\sim z$.
\item
$\lp \BP^{+}_{z}, \BP^-_{z'}\rp = 0$ unless $z \sim z'$.
\end{enumerate}
\item
The functions $\{\BQ_{z}^{\pm}\}$ are characterized as the dual
basis of $\{\BP_{z}^{\pm}\}$, i.e., 
\begin{equation*}
\lp \BP_{z}^+, \BQ_{z'}^-\rp = \lp \BQ_{z}^+, \BP_{z'}^-\rp
                               = \d_{z,z'} 
\qquad (z, z' \in \wt Z_{n,q}^{r,0}).
\end{equation*}
\end{enumerate}
\end{prop}
\begin{proof}
By substituting $t = 0$ into (3.4.2), we have
\begin{equation*}
\BOm_q(x,y;0) = \sum_{(\Ba,b)}z_{\Ba,b}\iv \Bp_{\Ba,b}(x)\ol{\Bp}_{\Ba,b}(y),
\end{equation*}
where $z_{\Ba,b} = |Z_W(w)|$ for $w = w_{\Ba}(b)$.
On the other hand, starting from the expansion of $\Om(x,y;0)$ by
means of Schur functions $s_{\Ba}(x)$ as given in Remark 4.9 in [S],
by using a similar argument as in the proof of (3.4.1), one can show
that
\begin{equation*}
\BOm_q(x,y;0) = \sum_{(\Ba,\f)}\Bs_{\Ba,\f}(x)\Bs_{\Ba,\f}(y).
\end{equation*}
It follows that one can define a hermitian form on $\Xi^n_{\CC}(p,q)$
satisfying the properties that $\lp \Bp_{\xi}, \Bp_{\xi'}\rp = \d_{\xi,\xi'}$
($\xi,\xi' \in \wt\CP_q^{n,e,p}$), and that  
$\lp\Bs_{z}, \Bs_{z'}\rp = \d_{z,z'}$ ($z, z' \in \wt\CP^q_{n,e,p}$). 
Hence the seaquiliner form on $\Xi^n_{\CC(t)}(p,q)$ is reduced to the
hermitian form on $\Xi^n_{\CC}(p,q)$ by substituting $t = 0$.
Now the arguments in the proof of Proposition 4.8 and Remark 4.9 in
[S] can be applied to our situation, and one can show that there
exist unique functions 
satisfying the properties (i), (ii) in the proposition, 
respectively.
So, we have only to show that 
$\BP^{\pm}_z, \BQ^{\pm}_z$ satisfy the properties in the
proposition.
\par
First we note that 
$\lp P^{+,j}_{\BLa,\f}, P^{-,j}_{\BLa',\f'}\rp = 0$ unless 
$\BLa \sim \BLa'$ by Proposition 3.7. It follows from (3.8.4)
that we have 
$\lp \BP^+_z,\BP^-_{z'}\rp = 0$ unless $z \sim z'$.  This shows
the proerty (b) of (i). Next we shall show the property (a).
Clearly the statement (a) is equivalent to the statement that  
\begin{equation*}
\tag{3.10.1}
\BP^{\pm}_{z} = \Bs_{z} + \sum_{z'}d_{z,z'}^{\pm}\BP^{\pm}_{z'}
\end{equation*}
for any $z,z' \in \wt Z_{n,q}^{r,0}$, 
where $d_{z,z'}^{\pm} = 0$ unless $z' \prec z$ and $z\not\sim z'$. 
\par
We consider the equation (3.10.1) for a fixed $z$ 
with respect to the $+$ sign. 
Take $z'$ such that 
$z' \prec z$ and $z \not\sim z'$.  Taking the scalar product with
$\BP^-_{z'}$ on both sides of (3.10.1), we have
\begin{equation*}
\tag{3.10.2}
\lp \Bs_z, \BP^-_{z'}\rp + 
       \sum_{z'' \sim z'}d^+_{z,z''}\lp \BP^+_{z''},\BP^-_{z'}\rp = 0.
\end{equation*}
When $z'$ runs over all the elements in a fixed similarity class, 
(3.10.2) can be regarded as a system of equations with unknown
variables $\{ d^+_{z, z''}\}$.  As in the arguments in the first part
of the proof (cf. Remark 4.9 in [S]), the matrix 
$(\lp \BP^+_{z''}, \BP^-_{z'}\rp)$ is non-singular.  Hence (3.10.2)
has a unique solution $\{ d^+_{z, z''}\}$.
Put
\begin{equation*}
F_z = \Bs_z + \sum_{z''}d^+_{z,z''}\BP^+_{z''}
\end{equation*} 
using thus determined $d^+_{z,z''}$ for $z'' \prec z$ and $z \not\sim z''$. 
Then clearly we have 
\begin{equation*}
\tag{3.10.3}
\lp \BP^+_z, \BP^-_{z'}\rp = \lp F_z, \BP^-_{z'}\rp
\qquad\text{ for  $z' \prec z$ and $z' \not\sim z$}. 
\end{equation*}
On the other hand, it follows from Proposition 3.7 that we have
\begin{equation*}
\lp P^+_{\BLa(\Ba)}, P^-_{\BLa'}\rp = \lp \Bs_{\Ba}, P^-_{\BLa'}\rp
\qquad \text{ for $\BLa(\Ba) \prec \BLa'$ or $\BLa(\Ba) \sim \BLa'$}.
\end{equation*}
This implies, in view of (3.8.4), that
\begin{equation*}
\tag{3.10.4}
\lp\BP^+_z, \BP^-_{z'}\rp = \lp\Bs_{z}, \BP^-_{z'}\rp 
= \lp F_z, \BP^-_{z'}\rp
\qquad\text{ for $z \prec z'$ or $z \sim z'$}.
\end{equation*} 
Now (3.10.3) and (3.10.4) implies that $\BP^+_z = F_z$, and (3.10.1)
is proved for $\BP^+_z$.  The proof for $\BP^-_z$ is similar.
\par
Finally, we show (ii).  By Proposition 3.7 (ii), we have 
\begin{equation*}
\lp P^{+,j}_z, Q^{-,j}_{z'}\rp_j = 
               \begin{cases}
     \f(\t^j)\ol{\f'(\t^j)}c_{\Ba,p}
          &\quad\text{ if $c_{\Ba,p} \mid j$ and $z = z'$}, \\
      0    &\quad\text{ otherwise}.
               \end{cases} 
\end{equation*}
This implies, by (3.8.4), that 
$\lp\BP^+_z, \BQ^-_{z'}\rp = \d_{z,z'}$.
The formula $\lp\BQ^+_z, \BP^-_{z'}\rp = \d_{z,z'}$ is shown 
similarly.
\end{proof}
\par
As a corollary to Proposition 3.10, we have
\begin{cor} %%%% Cor. 3.11 
$\BOm_q(x,y;t)$ has the following expansions.
\begin{align*}
\tag{3.11.1}
\BOm_q(x.y;t) &= \sum_{z,z'}
         b_{z, z'}(t)
           \BP_{z}^+(x;t)\bar \BP_{z'}^-(y;t)  \\
\tag{3.11.2}
\BOm_q(x,y;t) &= \sum_{z}
                 \BQ^+_{z}(x;t)\bar\BP^-_{z}(y;t)
     = \sum_{z}\BP^+_{z}(x;t)\bar\BQ^-_{z}(y;t),
\end{align*}
where $z = (\BLa,\f), z' = (\BLa',\f')$ run over all the elements
in $\bigcup_{n=1}^{\infty}\wt Z_{n,q}^{r,0}$, and 
$b_{z, z'} = 0$ unless $|\BLa| = |\BLa'|$ and
$\BLa \sim \BLa'$ in (3.9.1).  In (3.9.2), $z = (\BLa,\f)$ runs over 
all the elements in $\bigcup_{n=1}^{\infty}\wt Z_{n,q}^{r,0}$.
\end{cor}
\begin{proof}
(3.11.2) follows from (3.4.1) and Proposition 3.10 (ii) by 
a standard argument, (e.g. [M, I, 4.6]).  We show (3.11.1).
We write $\BP^-_z = \sum_{z'}a_{z,z'}\BQ^-_{z'}$.  
Since $\lp\BP^+_{z'}, \BP^-_z\rp = a_{z,z'}$, Proposition 3.10 (b) 
implies that  
the transition matrix between two bases $\{ \BP^-_z\}$ and
$\{ \BQ^-_z\}$ is a block diagonal matrix with respect to the 
partition by similarity classes.
Hence its inverse matrix is also block diagonal, i.e., 
$\BQ^-_z$ is written as $\BQ^-_z = \sum_{z'}b_{z,z'}\BP^-_{z'}$, 
where $b_{z,z'} = 0$ unless $z \sim z'$. 
Substituting this into the second formula in (3.11.2), we obtain 
(3.11.1).  Thus the corollary is proved.
\end{proof}
\par\bigskip
\section{Green functions associated to $G(e,p,n)$}
\para{4.1.}
For functions $f, h$ on $\s^q W$, we define an inner product
$\lp f,h\rp_q $ by
\begin{equation*}
\lp f, h\rp_q = |W|\iv\sum_{w \in \s^q W}f(w)\ol{h(w)}.
\end{equation*}
It is known in general that the number of $\s^q$-stable irreducible 
characters 
is equal to the number of $W$-orbits in $\s^q W$. An $\s^q$-stable
irreducible character $\x$ on $W$ can be extended to an irreducible 
character of $\lp\s^q\rp\ltimes W$, in $e/q$ distinct way.  
We fix an extension $\wt\x$ of $\x$.  Then $\wt\x$ gives a function 
on $\s^q W$ by restriction, which is constant on each $W$-orbit.  
It is known that
\begin{equation*}
\tag{4.1.1}
\lp \wt\x, \wt\x' \rp_{\, q} = \begin{cases}
                1 &\quad\text{ if } \wt\x = \wt\x', \\
                0 &\quad\text{ if }  \x \ne \x'.
                              \end{cases}
\end{equation*}
Now the set of $\s^q$-stable irreducible characters is parametrized 
by $\wt\CP^q_W$.  Also the set of $W$-orbits in $\s^q W$ is
parametrized by $\wt\CP_q^W$.  We fix a total order on 
$\wt\CP_W^q \simeq \wt Z_{n,q}^{0,0}$ as in 3.6.  Let  
$X = (\wt\x^z(\xi))\, (z \in \wt\CP_W^q, \xi \in \wt\CP_q^W)$ 
be the character table of $\s^q W$.
We define $H$ as the diagonal matrix indexed by $\wt\CP^q_{n,e,p}$
whose diagonal entry 
corresponding to $(\xi,\xi)$ is $z_{\xi}\iv = |Z_W(w_{\Ba}(b)|\iv$ 
for $\xi = (\Ba,b)$.   
Then (4.1.1) can be written in a matrix form
\begin{equation*}
\tag{4.1.2}
{}^tXH\bar X = I. 
\end{equation*}
\para{4.2.}
Let $M = \CC^n$ be the natural reflection representation of 
$\wt W$.  We consider $M$ as a $W$-module by the restrition.
Let $S(M)$ be the symmetric algebra of $M$, and $I_+$ the ideal
of $S(M)$ generated by $W$-invariant homogeneous vectors of 
strictly positiv degree.  We deonte by $R = S(M)/I_+$ the coinvariant
algebra of $W$. The Poincar\'e polynomial $P_W(t)$ associated to $W$
is defined in terms of the graded algebra $R = \bigoplus R_i$.  In our
case, $P_W(t)$ is explicitly given as 
\begin{equation*}
P_W(t) = \prod_{i=1}^{n-1}\frac{t^{ei}-1}{t-1}
               \cdot\frac{t^{dn}-1}{t-1}. 
\end{equation*}
Since $\s$ normalizes $W$, $\s$ acts naturally on $R$,
stablizing each summand $R_i$.  
We denote by $\wt R_i$ the thus obtained $\wt W$-module. 
\par
Let $\det_M$ be the linear character of $\wt W$ defined by the 
determinant on $M$. For each $\wt W$-stable function $f$ on $\s^qW$, 
we deifne $R_q(f) \in \CC(t)$ by
\begin{equation*}
\tag{4.2.1}
R_q(f) = (\z^{qd}t^{dn}-1)\prod_{i=1}^{n-1}(t^{ei}-1)
           \cdot\frac{1}{|W|}\sum_{w \in \s^qW}
              \frac{\det_M(w)f(w)}{\det_M(t\cdot\id_M - w)}.
\end{equation*}   
Then we have
\begin{equation*}
R_q(f) = \sum_i\lp f, \wt R_i\rp_{q}t^i,
\end{equation*} 
and, in particular, $R_q(\wt\x^z) \in \ZZ[\z][t]$ for 
$z \in \wt\CP_W^q$. 
\para{4.3.}
Let us define a square matrix $\Om_q = (\w_{z,z'})$ 
indexed by $\wt Z_{n,q}^{0,0}$ by
\begin{equation*}
\w_{z, z'} = t^{N^*}
            R_q(\wt\x^z\otimes\wt\x^{z}\otimes\ol{\det}_M),
\end{equation*}
where $N^*$ is the number of reflections in $W$.  
Also we define a matrix $\Om_q' = (\w'_{z, z'})$ by
\begin{equation*}
\tag{4.3.1}
\w'_{z, z'} = t^{N^*}
         R_q(\wt\x^{z}\otimes\ol{\wt\x^{z'}}\otimes\ol{\det}_M).
\end{equation*} 
In either of the formulas, $\ol{\det}_M$ etc. denote the complex
conjugates of them.  
Let $P = (p_{z, z'})$ and 
$\vL = (\la_{z, z'})$ be the matrices indexed by 
$\wt Z_{n,q}^{0,0}$.  We consider the following system of equations 
with unknown variables 
$p_{z, z'}, \la_{z, z'}$; 
\begin{align*}
\la_{z, z'} &= 0 \quad\text{ unless } z \sim z', \\
\tag{4.3.2}
p_{z, z'} &= 0   \quad\text{ unless either } z \succ z' 
                        \text{ and } z \not\sim z', 
                        \text{ or } z = z', \\
p_{z, z'} &= t^{a(z)}, \qquad\\
P\vL\,&{}^t\!P = \Om_q.
\end{align*}
As discussed in [S, 1.5], the system of equations (4.1.2) is
equivalent to the system of equations
\begin{equation*}
\tag{4.3.3}
P'\vL'\, {}^t\!P'' = \Om_q',
\end{equation*}
where either of $P' = (p'_{z, z'})$ and $P'' = (p_{z, z'})$ 
satisfies  similar conditions as in the second and third one in 
(4.3.2), and $\vL' = (\la_{z, z'})$ as in the first one.
(Actually, it is shown that $P = P'$).
In the reamider of this section, we shall show that the solution 
of (4.1.3) can be descrbed in a combinatorial way.
\par\medskip\noindent
\addtocounter{thm}{1}
{\bf Remark 4.4.} 
In the case where $e = p = 2$, $G(e,e,n)$ is the Weyl
group of type $D_n$. In this case, the equation (4.2.2) coincides with
(4.2.3), and the solution for it describes the Green functions of
$SO_{2n}$ defined over a finite field of odd characteristic
 of split type (resp. non-split type) if $q = 0$ (resp. $q=1$).
\para{4.5.}
For two bases $X, Y$ of $\Xi^n_{\CC(t)}(p,q)$, we denote by
$M(X,Y)$ the transition matrix between $X$ and $Y$.
Let $X_{\pm}(t)$ be the transition matrix $M(\Bp, \BP^{\pm})$ between
the power sum symmetric functions $\{ \Bp_{\xi} \}$ and the Hall-Littlewood
functions $\{ \BP_{z}^{\pm}\}$, i.e.,
\begin{equation*}
\tag{4.5.1}
\Bp_{\xi}(x) = \sum_{z}X^{z}_{\xi,\pm}(t)
                   \BP^{\pm}_{z}(x;t).
\end{equation*}
(Here $\xi \wt\CP_q^{n,e,p}$ is the row-index and 
$z = (\BLa,\f) \in \wt\CP^q_{n,e,p}$ is the column-index.) 
Then by (3.8.1) and by the properties of $P_{\BLa}$ shown in
[S], we see that  $X_{\xi;\pm}^{z}(t) \in \CC(t)$ and 
is equal to zero unless $|z| = |\xi|$.  
(For $z = (\Ba,\f)$ and $\xi = (\Bb,b)$, we put
$|z| = |\Ba|$ and $|\xi| = |\Bb|$, respectively.)
Moreover since 
$\BP^{\pm}_{z}(x;0) = \Bs_{z}(x)$ 
by (3.9.2), we have 
\begin{equation*}
\tag{4.5.2}
X_{\xi;\pm}^{z}(0) = \wt\x^{\Ba,\f}(w_{\Bb}(b))
\end{equation*} 
by Proposition 2.11.  It follows that  $X_{\pm}(0) = M(\Bp,\Bs)$ 
is the character table of $W$.  In particular, $X_{\pm}(0)$ is 
independent of the sign, which we  denote simply as $X(0)$. 
By combining (3.4.2) and (3.11.1), we have
\begin{equation*}
\sum_{|\xi| = n}z_{\xi}(t)\iv \Bp_{\xi}(x)\ol{\Bp}_{\xi}(y)
  = \sum_{z,z'} b_{z,z'}(t)
         \BP_z^+(x;t)\bar\BP_{z'}^-(y;t),
\end{equation*}
where in the right hand side, $z, z'$ 
run over all elements in $\wt Z_{n,q}^{r,0}$. 
We put $D(t) = D_- = (b_{z,z'}(t))$, the matrix indexed by 
$\wt\CP^q_{n,e,p} \simeq \wt Z_{n,q}^{0,0}$, and denote by $Z(t)$ the
diagonal matrix indexed by
$\wt\CP_q^{n,e,p}$, where the entry corresponding to $(\xi, \xi)$ 
is given by $z_{\xi}(t)$.   Substituting
(4.5.1) into the above equation, we have
\begin{equation*}
{}^tX_+(t)Z(t)\iv \bar X_-(t) = D(t).
\end{equation*}
If we put $\vL(t) = D(t)\iv$, the above 
formula is equivalent to
\begin{equation*}
\tag{4.5.3}
\bar X_-(t)\vL(t)\, {}^t\!X_+(t) = Z(t).
\end{equation*}
\para{4.6.}
Let $K_{\pm}(t) = M(\Bs, \BP^{\pm})$ be the transition matrix between
Schur functions and Hall-Littlewood functions, i.e., for 
$z, z' \in \wt\CP^q_{n,e,p}$, 
\begin{equation*}
\tag{4.6.1}
\Bs_{z'}(x) = \sum_{z}K^{\pm}_{z', z}(t)\BP^{\pm}_{z}(x;t).
\end{equation*}  
Then by Proposition 3.10, $K_{\pm}(t)$ is a block lower triagnlar matrix
with identity diagonal blocks, with entries in $\CC(t)$.
Since 
$M(\Bs,\BP^{\pm}) = M(\Bp,\Bs)\iv M(\Bp,\BP^{\pm})$, we have 
$K_{\pm}(t) = X(0)\iv X_{\pm}(t)$.  Substituting this into
(4.5.3), we have
\begin{equation*}
\tag{4.6.2}
K_-(t)\vL(t)\,{}^t\!K_+(t) = 
                \bar X(0)\iv Z(t)\,{}^t\!X(0)\iv.
\end{equation*} 
\par
We now define Green functions $Q^{z}_{\xi;\pm}(t) \in \CC(t)$
by
\begin{equation*}
Q^{z}_{\xi;\pm}(t) = t^{a(z)}X^{z}_{\xi;\pm}(t\iv).
\end{equation*}
If we put 
$\wt K_{z', z}^{\pm}(t) =
          t^{a(z)}K_{z', z}^{\pm}(t\iv)$,
$Q^{z}_{\xi;\pm}(t)$ can be written as
\begin{equation*}
Q^{z}_{\xi;\pm}(t) = \sum_{z''}
         \wt\x^{z''}(w_{\Bb}(b))
             \wt K^{\pm}_{z'', z}(t)
\end{equation*} 
for $\xi = (\Bb,b)$.
Let $\wt K_{\pm}(t) = (\wt K^{\pm}_{z', z}(t))$.  Then 
$\wt K_{\pm}(t) = K_{\pm}(t\iv)T$, where $T$ is a a diagonal matrix 
with diagonal entries 
$t^{a(z)}$.  Hence (4.6.2) can be rewritten as
\begin{equation*}
\tag{4.6.3}
\wt K_-(t)\vL'(t)\,{}^t\!\wt K_+(t) = 
           \bar X(0)\iv Z(t\iv)\,{}^t\!X(0)\iv,
\end{equation*}
where $\vL'(t) = T\iv\vL(t\iv)T$.
Here $\vL'(t)$ is still a block diagonal matrix,
and $\wt K_{\pm}(t)$ are block lower triangular matrices, 
where the diagonal blocks consist of scalar matrices $t^{a(z)}I$.
\par
Put $\wt\vL(t) = t^{-n}\GG(t)\vL'(t)$ with
\begin{equation*}
\GG(t) = t^{N^*}(\z^{qd}t^{dn}-1)\prod_{i=1}^{n-1}(t^{ei}-1).
\end{equation*}
The following result gives a combinatorial description 
of the solution of (4.3.3).
\par\medskip
\begin{thm}%%%Theorem 4.7
We have
\begin{equation*}
\wt K_-(t)\wt\vL(t)\,{}^t\!\wt K_+(t) = \Om_q'.
\end{equation*}
Hence $P' = \wt K_-(t), P'' = \wt K_+(t)$ and 
$\vL' = \wt\vL(t)$ give a solution for the
equation (4.3.3). 
\end{thm}
\begin{proof}
Let $Y$ be the right hand side of (4.6.2).  We shall compute 
$Y$.  
Let $w_{\Ba}(b)$ be the
element in $\s^qW$ corresponding to $\xi = (\Ba,b)$.
Since $X(0)$ is the
character table of $\s^qW$, we have ${}^t\!X(0)H\bar X(0) = I$ by 
(4.1.2), where $H$ is the diagonal matrix with diagonal entries 
$z_{\xi}\iv$.  It follows that
\begin{equation*}
\tag{4.7.1}
Y = {}^t\!X(0)H Z(t\iv)H\bar X(0).
\end{equation*}
Now it is known that if 
$\Ba = (\a^{(k)}_j) \in \CP_{n,e}$,  we have
\begin{equation*}
{\det}_M (t\cdot\id_M - w_{\Ba}(b)) = 
      \prod_{k = 0}^{e-1}\prod_{j=1}^{l(\a^{(k)})}(t^{\a^{(k)}_j} - \z^k).
\end{equation*}
Hence by (3.3.8), 
$z_{\Ba,b}(t\iv) = z_{\Ba,b}t^n{\det}_V(t\cdot\id_V - w_{\Ba}(b))\iv$.
In particular,  the $(z, z')$-entry of $Y$ is equal to
\begin{equation*}
t^n|W|\iv\sum_{w\in \s^qW}
 \frac{\wt\x^{z}(w)\ol{\wt\x^{z'}}(w)}
         {{\det}_M(t\cdot\id_M - w)}.
\end{equation*}
Therefore $\Om_q' = t^{-n}\GG(t)Y$ and the theorem follows.
\end{proof}
\para{4.8.}
Recall that Kostka functions $K^{\pm}_{\Ba,\Bb}(t)$ associated to symbols in
$Z_n^{r,0}(\Bm)$  is defined in [S] by the formula
\begin{equation*}
\tag{4.8.1}
s_{\Ba}(x) = \sum_{\Bb \in Z_n^{0,0}}
               K^{\pm}_{\Ba,\Bb}(t)P^{\pm}_{\BLa(\Bb)}(x;t).
\end{equation*}
By making use of the preceding discussions, we can describe the
Kostka functions $K^{\pm}_{z,z'}$ associated to $\wt Z_{n,q}^{r,0}$
in terms of $K^{\pm}_{\Ba,\Bb}$.  
For  $\Ba,\Ba' \in Z_n^{0,0}$, put $c = c_{\Ba,p}$, $c' = c_{\Ba',p}$.
We assume that $cq\equiv 0, c'q \equiv 0 \pmod p$.
For each $j$ such that $c' \mid j$, put
\begin{equation*}
\tag{4.8.2}
L^{j,\pm}_{\Ba,\Ba'} = \begin{cases}
       \sum_{i=0}^{c'-1}\z^{-qid}K^{\pm}_{\Ba\{ j\},\th^i(\Ba')\{ j\}}
                         &\quad\text{ if } c \mid j, \\
      0                  &\quad\text{ otherwise}.
                       \end{cases}
\end{equation*}
We denote by $L_{\Ba,\Ba'}^{\pm}$ the diagonal matrix of size
$|\vG_{\Ba'}|$ whose $ii$-entry is given by $L^{c'i,\pm}_{\Ba,\Ba'}$.
Let $\F' = (\f'(\t^j))_{\f',j}$ be the character table of
$\vG_{\Ba'}$ ($\f' \in \vG\wg_{\Ba'},  c' | j)$, and let
$\F = (\f(\t^j))_{\f, j}$ be a part of the character table of
$\vG_{\Ba}$, i.e., $\f \in \vG\wg_{\Ba}$, but $j$ is taken under the
condition that $c\mid j$ and $c'\mid j$.  For fixed $\Ba, \Ba'$ as above,
we consider the matrix $(K^{\pm}_{\f,\f'})$, indexed by 
$\f \in \vG\wg_{\Ba}, \f' \in \vG\wg_{\Ba'}$, defined by 
$K^{\pm}_{\f,\f'} = K^{\pm}_{z,,z'}$. 
We have the following prpopsition.
%%%
\begin{prop}%%%%Prop.4.9
The matrix $(K^{\pm}_{\f,\f'})$ can be written as 
\begin{equation*}
(K^{\pm}_{\f,\f'}) = |\vG_{\Ba'}||\vG_{\Ba}|\iv\F L^{\pm}_{\Ba,\Ba'}{\F'}\iv.
\end{equation*}
\end{prop}
\begin{proof}
By using the relation (4.8.1) for Schur functions with respect to the
variables $\CX_j$, we can exprese $s_{\Ba\{ j\}}$ in terms of a linear
combination of $P^{\pm}_{\Bb}$ with $\Bb \in Z_{n/h_j}^{0,0}$.  
Then by applying $\vT_q$ on
both sides, we have
\begin{equation*}
\vT_q(s_{\Ba\{ j\}}) = \sum_{\Ba'}K^{\pm}_{\Ba\{ j\},\Ba'\{ j\}}
   \vT_q(P^{\pm}_{\BLa(\Ba')\{ j\}}), 
\end{equation*}
where $c' = c_{\Ba',p}$, and $\Ba'$ in the right hand side runs 
over all $\Ba' \in Z_n^{0,0}$ such that $j \mid c'$  and that 
$c' \equiv 0 \pmod p$.  Hence if we choose a representative
$\Ba'$ from each $\th$-orbit, the above formula turns out to be
\begin{equation*}
\tag{4.9.1}
\vT_q(s_{\Ba\{ j\}}) = \sum_{\Ba' \in Z_n^{0,0}\ssim_p}
       \sum_{i=0}^{c'-1}\z^{-qid}K^{\pm}_{\Ba\{ j\},\th^i(\Ba')\{ j\}}
          \vT_q(P^{\pm}_{\BLa(\Ba')\{ j\}}).
\end{equation*}
On the other hand, by using (4.6.1), we can express $\Bs_z$ in terms
of a linear combination of $\BP_{z'}^{\pm}$.  In particular,  we have
\begin{equation*}
s_z^j = \sum_{z'}K^{\pm}_{z,z'}P^{\pm,j}_{z'}
\end{equation*}
fpr any $j$ such that $0 \le j < p$.
Now by making use of (3.3.5) and a similar formula for $P_{\BLa}$, 
we can rewrite the above formula as
\begin{equation*}
|\vG_{\Ba}|\iv\f(\t^j)\vT_q(s_{\Ba\{j\}}) = 
\sum_{\Ba' \in Z_n^{0,0}\ssim_p}|\vG_{\Ba'}|\iv
     \sum_{\f'\in \vG\wg_{\Ba'} }K^{\pm}_{z,z'}\f'(\t^j)
               \vT_q(P^{\pm}_{\BLa(\Ba')\{ j\}}) 
\end{equation*}
for any $j$ such that $c'\mid j$.  Note that the left hand side 
is understood to be 0 in the case where $j$ is not a multiple of $c$. 
Since $\vT_q(P^{\pm}_{\BLa(\Ba'\{j\}})$ are linearly independent in 
$\Xi_{p,q}^{n/h_j}(\CX_j)$, by comparing with (4.9.1), we have
\begin{equation*}
\tag{4.9.2}
|\vG_{\Ba}|\iv\f(\t^j)\sum_{i=0}^{c'-1}\z^{-qid}
              K^{\pm}_{\Ba\{j\}, \th^i(\Ba')\{j\}}
  = |\vG_{\Ba'}|\iv\sum_{\f'}K_{\f,\f'}^{\pm}\f'(\t^j)
\end{equation*}  
for each $\Ba'$. Now, (4.9.2) holds for any $\f \in \vG_{\Ba}\wg$ and
any $j$ such that $c' \mid j$.  Hence it can be translated to a matrix
equation
\begin{equation*}
|\vG_{\Ba}|\iv\F L_{\Ba,\Ba'}^{\pm} = |\vG_{\Ba'}|\iv(K^{\pm}_{\f,\f'})\F'. 
\end{equation*}  
This proves the proposition.
\end{proof}   
The following special case would be worth mentioning.
\begin{cor}%%%Cor 4.10
Assume that $\vG_{\Ba} = \vG_{\Ba'} = \{1\}$.  Then we have
\begin{equation*}
K^{\pm}_{(\Ba,1), (\Ba',1)} = \sum_{i=0}^{p-1}\z^{-qid}
             K^{\pm}_{\Ba, \th^i(\Ba')}.
\end{equation*}
\end{cor}
%%%%
As in the case of Green functions associated to $G(e,1,n)$, one can
expect that, for $z,z' \in Z_{n,q}^{r,0}$, $K^{\pm}_{z,z'}(t)$ is 
a polynomial in $t$ with positive integral coefficeints.  In view of 
Corollary 4.10, the conjecture for $G(e,1,n)$, i.e., 
$K^{\pm}_{\Ba, \Ba'}(t) \in \ZZ[t]$,  implies that 
$K^{\pm}_{z,z'}(t) \in K[t]$ with $K = \QQ(\z)$.   
\par\medskip\noindent
{\bf Remarks 4.11.} (i) \ Proposition 4.9 asserts that Green functions
associated to $G(e,p,n)$ are completely described in terms of Green
functions associated to $G(e',1,n')$ with various 
$e'\mid e, n'\mid n$.  
In the classical situation, 
this means that Green functions of type $D_n$ can be described 
in terms of Green functions assocaited to Weyl groups of type $B_n$,
and $\FS_{n/2}$ for even $n$.
However, note that this does not mean such Green functions are
described in terms of Green functions of $Sp_{2n}$ or $SO_{2n+1}$. 
In fact, the Green functions of $Sp_{2n}$ or $SO_{2n+1}$ 
are defined in terms of symbols of the
type $Z_n^{r,s}(\Bm)$ with $\Bm = (m+1, m)$, while our Green functions
are related to the symbols of the type $Z_n^{r,0}(\Bm)$ with 
$\Bm = (m,m)$, and the structure of similarity classes is different.
\par
(ii) \ 
The group $G(e,e,2)$ coinicdes with the dihedral group of degree $2e$.
In particular, it coincides with the Weyl group of type $A_2, B_2$ or
$G_2$ for $e = 3,4$ or 6, respectively.  The case where $e=3$, our
Green function coincides with the Green function asociated to 
$GL_3(\Fq)$.   
However, as Lemma 5.4 in the next section shows,  the Green function 
associated to $G(e,e,2)$ for $e =4$ or $e=6$ (for any $r$) does not 
coincide with the Green function of type $B_2$ or $G_2$. 
\par\bigskip
\par\medskip
\section{Examples}
\para{5.1.}
In this section, we give some examples of the matrices 
$P = \wt K_-(t)$ associated to the group $G(e,e,n)$.  
But before doing it, we shall give a general remakr in connection with
fake degrees.  Assume that $W = G(e,p,n)$, and we consider the case 
where $q = 0$.  We write $P = (p_{z,z'})$ for 
$z,z' \in \wt\CP_W^0$.
For 
$j = 0, \dots, e-1$, put $\Bb_j = (\b^{(0)}, \dots, \b^{(e-1)})$, 
where $\b^{(j)} = (1^n)$ and $\b^i = \emptyset$ for $i \ne j$. 
In this case, $\vG_{\Bb_j} = \{1\}$, and we identify $z = (\Bb_j,1)$
with $\Bb_j$.
Then among $\wt\CP_W^0$, 
$\CF = \{ \BLa(\Bb_0), \dots, \BLa(\Bb_{d-1})\}$ gives rise to a
similarity class in $\wt Z_{n,0}^{r,0}$, which has the maximum 
$a$-value.
\par
For each $z \in \wt\CP_W^0$, $R_0(\x^{z})$ coincides with the fake 
degree $\sum_{i \ge 0}\lp \x^{z}, R_i\rp t^i$ (cf. 4.2).  
The following formula shows that the fake degeree can be read from the 
matrix $P$, which is an analogue of Lemma 7.2 in [S].  The proof is 
completely similar to it.
\begin{lem} %%% lemma 5.2.
Let $W = G(e,p,n)$ and assume that $q=0$. Put $b = a(\Bb_0)$.  Then we have
\begin{equation*}
\sum_{j =0}^{d-1}t^{-b}R(\x^{\Bb_j})p_{z,\Bb_j} = R(\x^{z}),
\end{equation*}
where $R(\x^{\Bb_j}) = t^{n(n-1)/2 + jn}$.
\par
In particular, in the case where $W = G(e,e,n)$, we have 
$\CF = \{ \BLa(\Bb_0)\}$, and 
$p_{z, \Bb_0} = t^{n(n-1)/2 -b}R(\x^{z})$. 
\end{lem}  
\para{5.3.}
Assume that $W = G(e,e,2)$, the dihedral group of order $2e$.
We shall determine the matrix $P$ in the case where $q=0$.
Let $\Ba_0, \Ba_{ij}$ be $e$-partitions defined in 2.16.
The set $\wt\CP_W^0$ is given as follows.  
\begin{equation*}
\Ba_{0,0}, \quad \Ba_{0,j} \quad  (1 \le j < e/2), \quad
\Ba_{0,e/2}, \quad \Ba'_{0,e/2}, \quad \Ba_0, 
\end{equation*}
where in the case of $\vG_{\Ba} = \{ 1\}$, we write as $\Ba$ 
to denote the element $(\Ba, 1) \in \wt\CP_W^0$.  Moreover, we have
$\vG_{\Ba} = \ZZ/2\ZZ$ for $\Ba = \Ba_{0,e/2}$ (the case where
$e$ is even), and in this case,
we write as $\Ba_{0,e/2}, \Ba'_{0,e/2}$ instead of  
$(\Ba_{0,e/2},1), (\Ba_{0.e/2}, -1)$, respectively. 
The similarity classes in 
$\wt Z_{n,0}^{r,0}$ are given as follows.
\begin{equation*}
\CF_1 = \{ \BLa(\Ba_{0,0})\}, \quad  
\CF_2 = \{ \BLa(\Ba_{0,j}) \mid 1 \le j \le e/2\}, \quad
\CF_3 = \{ \BLa(\Ba_0)\},  
\end{equation*}
where in $\CF_2$, $\BLa(\Ba_{0,j})$ for $j = e/2$
is understood as 
$\BLa(\Ba_{0,e/2}), \BLa(\Ba'_{0,e/2})$. The $a$-values on 
similarity classes are given as follows. 
\begin{equation*} 
a(\CF_1) = e, \quad a(\CF_2) = 1, \quad a(\CF_3) = 0.
\end{equation*}
Note that the similarity classes and $a$-values of them are
independent from the choice of $r$.
\par
Next we describe irreducible characters $\x^z$ associated to 
$z \in \wt\CP_W^0$.  If $z = \Ba_{0,0}, \Ba_0$, $\x^z$ coincide with
the sign character $\ve$ and the identity character $1$ of $W$, 
respectively.  If $e$ is even, $\Ba_{0,e/2}$ and $\Ba'_{0,e/2}$ 
correspond two linear characters (long sign character $\e$ 
and short sigh character $\e'$)
of $W$.  $z = \Ba_{0,1}$ corresponds to the reflection character
$\r_1$ of $W$ (hence $\deg \r_1 = 2)$.  Other characters 
corresponding to $\Ba_{0,j}$ have all degree 2.
We remark that the above partition of $W\wg$ 
into similarity classes together with $a$-values coincide with 
the partition of $W\wg$ into families and their $a$-values 
introduced by Lusztig [L1] in connection with the classification of 
unipotent characters of finite reductiove groups, 
where $1, \ve, \r_1$ are the special characters 
corresponding to families. 
\par
The fake degrees $R(\x^z)$  of $W$ are given as follows.
\begin{align*}
\tag{5.3.1}
R(1) &= 1, \quad R(\r_j) = t^j + t^{e-j} \, (1 \le j < e/2) \\
R(\e) &= t^{e/2}, \quad R(\e') = t^{e/2},  \\ 
R(\ve) &= t^e. 
\end{align*} 
\par
We now define an order on the set $W\wg$ by 
\begin{alignat*}{2}
&\ve, \r_1, \r_2, \dots, \r_{m}, 1 &\quad &\text{ if } e = 2m+1, \\
&\ve, \r_1, \r_2, \dots, \r_{m-1}, \e, \e', 1
                                   &\quad &\text{ if } e = 2m,
\end{alignat*}
and define the order on $\wt\CP_W^0$ accordingly.
We consider the matrix equation $P\vL\,{}^tP = \Om$ indexed by
$\wt\CP_W^0$. Following
the partition of $\wt\CP_W^0$ into similarity classes, we regard 
$P, \vL, \Om$ as block matrices.  By our assumption, 
$P$ and $\vL$ have the following shape.
\begin{equation*}
P = \begin{pmatrix}
        t^e & 0 & 0 \\
        P_{21} & tI_k & 0 \\
        p_{31} & P_{32} & 1
    \end{pmatrix},
\qquad
\vL = \begin{pmatrix}
         \la_{11} & 0 & 0 \\
         0  & \vL_{22} & 0 \\
         0 & 0 & \la_{33}
      \end{pmatrix}  
\end{equation*} 
where $k$ is equal to $e/2-1$ (resp. $e/2 +1$) if $e$ is odd (resp.
even), $\vL_{22}$ is a matrix of degree $k$.  $P_{21}, P_{32}$ are
matrices of the shape $(k,1), (1,k)$, respectively.
\par
Since the matrix equation $P\vL\,{}^tP = \Om$ is completely the same 
as the matrix equation arising from the partition into families
studied by Geck and Malle [GM], where the case of dihedral groups is
explcitly computed,  the matrix $P$ in our case is described as
follows.   Note that in view of Lemma 5.2, the first column
of $P$ coincides with the fake degrees of $W$.
%%%%
\begin{lem}[{[GM]}]%%%Lemma 5.4.
Let $W = G(e,e,2)$ be the dihedral group of degree $2e$.  Assume that
$q = 0$.  Then the matrix $P$ associated to 
$\wt Z_{n,0}^{r,0}$ is given as follows.
\begin{equation*}
{}^tP_{21} = \begin{cases} 
(t + t^{e-1}, t^2 + t^{e-2}, \dots, t^m+ t^{m+1})
                       &\quad\text{ if $e = 2m+1$}, \\
(t + t^{e-1}, t^2 + t^{e-2}, \dots, 
         t^{m-1}+ t^{m+1}, t^{m}, t^{m})
                        &\quad\text{ if $e = 2m$.}
              \end{cases}  
\end{equation*}
and
$p_{31} = 1, \quad P_{32} = (1, 0, \dots, 0)$.
\par
Moreover, a part of the matrix $\vL$ is given
as follows.
\begin{equation*}
\la_{11} = 1, \quad (\vL_{22})_1 = t\iv(t^e-1)\cdot{}^tP_{21},
\quad \la_{33} = t^{e-2}(t^2-1)(t^e-1), 
\end{equation*}
where $(\vL_{22})_1$ denotes the first row of the matrix $\vL_{22}$.
\end{lem}
\para{5.5.}   
In the remainder of this section, we shall give some more
examples of Green functions assciated to $G(e,e,n)$.
Throughout those examples, 
we always assume that $r = 2$ and $q = 0$.
\par
First assume that $W  = G(3,3,3)$.  Then $\wt\CP^0_W$ is given
as follows.
\begin{align*}
\Ba_1 &= (1^3;-;-), \Ba_2 = (1^2;1;-), \Ba_3 = (1;1^2;-), \\  
\Ba_4 &= (1;1;1), \Ba_4' = (1;1;1)', \Ba_4'' = (1;1;1)'', \\ 
\Ba_5 &= (21;-;-), \Ba_6 = (2;1;-), \Ba_7 = (1;2;-), 
\Ba_8 = (3;-;-).
\end{align*} 
Note taht $|\vG_{\Ba_4}| = 3$, and we write $\Ba_4, \Ba_4', \Ba_4''$
instead of $(\Ba,\f)$ with $\f \in \vG_{\Ba_4}$. $\vG_{\Ba} = \{1\}$
for other $\Ba$.  Then the symbols and
similarity classes in $\wt Z_{n,0}^{2,0}$ are given as
\begin{align*}
\CF_1 &= \{ \BLa_1 = (531;420;420)\}, \\
\CF_2 &= \{ \BLa_2 = (31;30;20), \BLa_3 = (30;31;20)\}, \\
\CF_3 &= \{ \BLa_4 = (1;1;1), \BLa_4' = (1;1;1)', \BLa_4'' = (1;1;1)''\}, \\
\CF_4 &= \{ \BLa_5 = (41;20;20)\}, \quad
    \CF_5 = \{ \BLa_6 = (2;1;0), \BLa_7 = (1;2;0)\}, \\ 
\CF_6 &= \{ \BLa_8 = (3;0;0)\}.
\end{align*}
The matrices $\vL(\CF_i)$ are given as follows;
\begin{align*}
\vL(\CF_1) &= (1), \qquad 
\vL(\CF_2) = (t^6-1)\begin{pmatrix}
                2t^3+1 & t(t^3+2) \\
                t(t^3+2) & 3t^2
             \end{pmatrix} \\
\vL(\CF_3) &= t^3(t^3-1)(t^6-1)I_3, \qquad 
\vL(\CF_4) =  (t^3(t^3-1)(t^6-1)), \\
\vL(\CF_5) &= t^3(t^3-1)^2(t^6-1)\begin{pmatrix}
                  2 & t \\
                  t & 0
              \end{pmatrix}, \\
\vL(\CF_6) &= (t^6(t^3-1)^2(t^6-1)),
\end{align*}
where $I_3$ denotes the identity matrix of degree 3.
\par
The matrix $P$ of Green functions is given in Table 1.
\begin{table}[h]
\caption{$G(3,3,3) \quad  (r = 2)$}
\label{D_3: D_3}
\begin{center}
\begin{tabular}{|c||c|cc|ccc|c|cc|c|}
\hline 
 $(1^3;-;-)$ &  $t^9$ & & & & & & & & &  \\
 \hline
 $(1^2;1;-)$ & $2t^7+t^4$ & $t^4$ &  & & & & & & &  \\
 $(1;1^2;-)$ & $t^8+2t^5$ &  & $t^4$ & & & & & & &    \\
\hline
 $(1;1;1)$ & $t^6+t^3$ & $t^3$ &  & $t^3$ & & & & & &\\
 $(1;1;1)'$ & $t^6+t^3$ & $t^3$ &  &  & $t^3$ & & & & &\\ 
 $(1;1;1)''$ & $t^6+t^3$ & $t^3$ &  &  &  & $t^3$ & & & &\\
\hline
$(21;-;-)$ & $t^6+t^3$ & $t^3$ &  &  &  &  & $t^3$ & & &\\
\hline
$(2;1;-)$ & $2t^4 +t$ & $t$ & $t^3$ & $t$ & $t$ & $t$ & $t$ & $t$ & & \\ 
$(1;2;-)$ & $t^5+2t^2$ & $2t^2$ &  & $t^2$ & $t^2$ & $t^2$ & $t^2$
          &  & $t$ &  \\
\hline
$(3;-;-)$ & 1 & 1 &  & 1 & 1 & 1 & 1 & 1 &  & 1 \\
\hline
\end{tabular}
\end{center}
\end{table}
\para{5.6.}
Assume that $W = G(4,4,3)$.  Then $\wt\CP^0_W$ is given as
\begin{align*}
\Ba_1 &= (1^3;-;-;-), \quad 
\Ba_2 = (1^2;1;-;-), \quad \Ba_3 = (1^2;-;1;-), \\
\Ba_4 &= (1;1^2;-;-), \quad \Ba_5 = (21;-;-;-), \quad \Ba_6 = (1;1;1;-), \\
\Ba_7 &= (2;1;-;-), \quad \Ba_8 = (2;-;1;-), \quad
\Ba_9 = (1;2;-;-), \\
\Ba_{10} &= (3;-;-;-).  
\end{align*}
The symbols and similarity classes in $\wt Z_{n,0}^{2,0}$ are given as
\begin{align*}
\CF_1 &= \{ \vL_1 = (531;420;420;420)\}, \\ 
\CF_2 &= \{ \vL_2 = (31;30;20;20), \vL_3 = (31;20;30;20), 
           \vL_4 = (30;31;20;20)\}, \\
\CF_3 &= \{ \vL_5 = (41;20;20;20)\}, \qquad
\CF_4 = \{ \vL_6 = (1;1;1;0)\}, \\
\CF_5 &= \{ \vL_7 = (2;1;0;0), \vL_8 = (2;0;1;0), \vL_9 = (1;2;0;0)\}, \\
\CF_6 &= \{ \vL_{10} = (3;0;0;0)\}.
\end{align*}
The matrices $\vL(\CF_i)$ are given as follows.
\begin{align*}
\vL(\CF_1) &= (1), \\ 
\vL(\CF_2) &= (t^8-1)\begin{pmatrix}
(t^5+t^4+1) & t(t^5+t+1) & t(t^4+t^2+1)  \\  
t(t^5+t+1)  & t^2(t^2+t+1)  & t^2(t^2+t+1)  \\
t(t^4+t^2+1)  &  t^2(t^2+t+1)  &  t^2(t^4+t^2+1)
             \end{pmatrix}, \\
\vL(\CF_3) &= (t^5(t^3-1)(t^8-1)), 
\quad
\vL(\CF_4) = (t^4(t^2+t+1)(t^4-1)(t^8-1)), \\ 
\vL(\CF_5) &= t^5(t^3-1)(t^4-1)(t^8-1)\begin{pmatrix}
t+1 & t^2 & t \\
t^2 & 0 & 0 \\
t & 0 & t^2
          \end{pmatrix}, \\
\vL(\CF_6) &= (t^9(t^3-1)(t^4-1)(t^8-1)).         
\end{align*}
The matrix of $P$ of Green functions is given in Table 2.
\begin{table}[h]
\caption{$G(4,4,3) \quad  (r = 2)$}
\label{D_3: D_3}
\begin{center}
\begin{tabular}{|c||c|ccc|c|c|ccc|c|}
\hline 
 $(1^3;-;-;-)$ &  $t^{12}$ & & & & & & & & &  \\
 \hline
 $(1^2;1;-;-)$ & $t^{10}+t^9 + t^5$ & $t^5$ &  & & & & & & &  \\
 $(1^2;1;-;1)$ & $t^{11}+t^7 + t^6$ &  & $t^5$ & & & & & & &    \\
 $(1;1^2;-;-)$ & $t^{10}+t^8 +t^6$ &  &  & $t^5$ & & & & & &\\
\hline
 $(21;-;-;-)$ & $t^8+t^4$ & $t^4$ &  &  & $t^4$ & & & & &\\ 
\hline
 $(1;1;1;-)$ & $t^9+t^8+t^7+t^5+t^4+t^3$ & $t^4+t^3$ &  &$t^4$  &  & $t^3$ & & & &\\
\hline
$(2;1;-;-)$ & $t^6+t^5$ & $t$ & $t^4$ &  & $t$ &$t$  & $t$ & & &\\
$(2;-;1;-)$ & $t^7 +t^3+t^2$ & $t^3+t^2$ &  &  & $t^3$ & $t^2$ & & $t$ & & \\ 
$(1;2;-;-)$ & $t^6+t^4+t^2$ & $t^2$ &  & $t^3$ & $t^2$ & $t^2$ & &  & $t$ &  \\
\hline
$(3;-;-;-)$ & 1 & 1 &  &  & 1 & 1 & 1 &  &  & 1 \\
\hline
\end{tabular}
\end{center}
\end{table}

\end{document}